\documentclass[aos,MSNbibl,nameyear,dvips]{arximspdf}

\doi{10.1214/12-AOS998} 
\volume{40}
\issue{5}
\pubyear{2012}
\firstpage{2389}
\lastpage{2420}

\makeatletter
\renewcommand{\overline}{\bar}

\newcommand{\rrVert}{\Vert}
\newcommand{\llVert}{\Vert}
\newtheorem{theorem}{Theorem}
\newtheorem{corollary}{Corollary}
\newtheorem{proposition}{Proposition}
\newtheorem{lemma}{Lemma}
\newproclaim{remark}{Remark}
\newcommand{\eqref}[1]{(\ref{#1})}
\makeatother

\begin{document}
\begin{frontmatter}

\title{Optimal rates of convergence for sparse covariance matrix estimation}
\runtitle{Estimating sparse covariance matrix}

\begin{aug}
\author[A]{\fnms{T. Tony} \snm{Cai}\thanksref{t1}\ead[label=e1]{tcai@wharton.upenn.edu}\ead[label=u1,url]{http://www-stat.wharton.upenn.edu/\textasciitilde tcai}} \and
\author[B]{\fnms{Harrison H.} \snm{Zhou}\corref{}\thanksref{t2}\ead[label=e2]{huibin.zhou@yale.edu}\ead[label=u2,url]{http://www.stat.yale.edu/\textasciitilde hz68}}
\thankstext{t1}{Supported in part by NSF FRG Grant DMS-08-54973.}
\thankstext{t2}{Supported in part by NSF Career Award DMS-06-45676 and
NSF FRG Grant DMS-08-54975.}
\runauthor{T. T. Cai and H. H. Zhou}
\affiliation{University of Pennsylvania and Yale University}
\address[A]{Department of Statistics\\
The Wharton School\\
University of Pennsylvania\\
Philadelphia, Pennsylvania 19104\\
USA\\
\printead{e1}\\
\printead{u1}}
\address[B]{Department of Statistics\\
Yale University\\
New Haven, Connecticut 06511\\
USA\\
\printead{e2}\\
\printead{u2}}
\end{aug}

\received{\smonth{12} \syear{2010}}
\revised{\smonth{3} \syear{2012}}

%
\begin{abstract}
This paper considers estimation of sparse covariance matrices and
establishes the
optimal rate of convergence under a range of matrix operator norm and
Bregman divergence losses. A major focus is on the
derivation of a rate sharp minimax lower bound. The problem exhibits new
features that are significantly different from those that occur in the
conventional nonparametric function estimation problems. Standard techniques
fail to yield good results, and new tools are thus needed.

We first develop a lower bound technique that is
particularly well suited for treating ``two-directional'' problems such as
estimating sparse covariance matrices. The result can be viewed as a
generalization of Le Cam's method in one direction and Assouad's Lemma in
another. This lower bound technique is of independent interest and can be
used for other matrix estimation problems.

We then establish a rate sharp minimax lower bound for estimating sparse
covariance matrices under the spectral norm by applying the general lower
bound technique. A thresholding estimator is shown to attain the optimal
rate of convergence under the spectral norm. The results are then
extended to the general matrix $\ell_w$ operator norms for $1\le w\le
\infty$. In addition, we give a unified result on the minimax rate of
convergence for sparse covariance matrix estimation under a class of
Bregman divergence losses.
\end{abstract}

%
%

\end{frontmatter}

\section{Introduction}

Minimax risk is one of the most widely used benchmarks for optimality, and
substantial efforts have been made on developing minimax theories in the
statistics literature.\vadjust{\goodbreak}
A key step in establishing a minimax theory is the derivation of minimax
lower bounds and several effective lower bound arguments based on hypothesis
testing have been introduced in the literature. Well-known techniques
include Le Cam's method, Assouad's lemma and Fano's lemma. See
\citet{LeCam86}
and \citet{Tsybakov} for more detailed discussions on minimax
lower bound
arguments.

Driven by a wide range of applications in high dimensional data analysis,
estimation of large covariance matrices has drawn considerable recent
attention. See, for example, \citeauthor{BicLev08N1} (\citeyear
{BicLev08N1,BicLev08N2}), \citet{ElK08},
\citet{Ravetal08}, \citet{LamFan09}, \citet
{CaiZhaZho10} and \citet{CaiLiu11}. Many
theoretical results, including consistency and rates of convergence, have
been obtained. However, the optimality question remains mostly open in the
context of covariance matrix estimation under the spectral norm, mainly due
to the technical difficulty in obtaining good minimax lower bounds.

In this paper we consider optimal estimation of sparse covariance matrices
and establish the minimax rate of convergence under a range of matrix
operator norm and Bregman divergence losses. A major focus is on the
derivation of a rate sharp
lower bound under the spectral norm loss. Conventional lower bound
techniques such as the ones mentioned earlier are designed and well suited
for problems with parameters that are scalar or vector-valued. They have
achieved great successes in solving many nonparametric function estimation
problems which can be treated exactly or approximately as estimation of
a finite or
infinite dimensional vector and can thus be viewed as
``one-directional'' in terms of the lower bound
arguments. In
contrast, the problem of estimating a sparse covariance matrix under the
spectral norm can be regarded as a truly
``two-directional'' problem where one direction is
along the
rows and another along the columns. It cannot be essentially reduced to a
problem of estimating a single or multiple vectors. As a consequence,
standard lower bound techniques fail to yield good results for this matrix
estimation problem. New and more general technical tools are thus needed.

In the present paper we first develop a minimax lower bound technique that
is particularly well suited for treating ``two-directional'' problems such
as estimating sparse covariance matrices. The result can be viewed as a
simultaneous generalization of Le Cam's method in one direction and
Assouad's lemma in another. This general technical tool is of independent
interest and is useful for solving other matrix estimation problems
such as
optimal estimation of sparse precision matrices.

We then consider specifically the problem of optimal estimation of sparse
covariance matrices under the spectral norm. Let $\mathbf
{X}_{1},\ldots,%
\mathbf{X}_{n}$ be a random sample from a $p$-variate distribution with
covariance matrix $\Sigma= ( \sigma_{ij} )_{1\leq i,j\leq p}$.
We wish to estimate the unknown matrix $\Sigma$ based on the sample $\{
\mathbf{X}_{1},\ldots,\mathbf{X}_{n}\}$. In this paper we shall use the
weak $%
\ell_{q}$ ball with $0\leq q<1$ to model the sparsity of the covariance
matrix $\Sigma$. The weak $\ell_{q}$ ball was originally used in
\citet{Abretal06} for a sparse normal means
problem. A weak $\ell_{q}$ ball of radius $c$ in $\mathbb{R}^{m}$ contains
elements with fast decaying ordered magnitudes of components,
\[
B_{q}^{m}(c)= \bigl\{ \xi\in\mathbb{R}^{m}\dvtx|
\xi|_{(k)}^{q}\leq ck^{-1}, \mbox{ for all $k=1, \ldots,
m$} \bigr\},
\]
where $|\xi|_{(k)}$ denotes the $k$th largest element in magnitude of the
vector $\xi$. For a covariance matrix $\Sigma=(\sigma_{ij})_{1\leq
i,j\leq p}$, denote by $\sigma_{-j,j}$ the $j$th column of $\Sigma$
with $%
\sigma_{jj}$ removed. We shall assume that $\sigma_{-j,j}$ is in a
weak $%
\ell_{q}$ ball for all $1\leq j\leq p$. More specifically, for $0\leq q<1$,
we define the parameter space $\mathcal{G}_{q}(c_{n,p})$ of covariance
matrices by
%
%
\begin{equation}
\mathcal{G}_{q}(c_{n,p})= \bigl\{ \Sigma= (
\sigma_{ij} )_{1\leq
i,j\leq p}\dvtx\sigma_{-j,j}\in
B_{q}^{p-1}(c_{n,p}), 1\leq j\leq p \bigr\} .
\label{sparseparaspace}
\end{equation}
In the special case of $q=0$, a matrix in $\mathcal{G}_{0}(c_{n,p})$
has at
most $c_{n,p}$ nonzero off-diagonal elements on each column.

The problem of estimating sparse covariance matrices under the spectral
norm has been
considered, for example, in \citet{ElK08}, \citet{BicLev08N2},
\citet{RotLevZhu09} and \citet{CaiLiu11}. Thresholding methods
were introduced, and rates of convergence in probability were obtained for
the thresholding estimators. The parameter space $\mathcal{G}_{q}(c_{n,p})$
given in (\ref{sparseparaspace}) also contains the uniformity class of
covariance matrices considered in \citet{BicLev08N2} as a special
case. We assume that the distribution of the $X_{i}$'s is subgaussian
in the
sense that there is $\tau>0$ such that
%
%
\begin{equation}
\mathbb{P}\bigl\{\bigl|\mathbf{v}^{T}(\mathbf{X}_{1}-\mathbb{E}
\mathbf{X}_{1})\bigr|>t\bigr\}\leq e^{-t^{2}/ ( 2\tau) }\qquad\mbox{for all
}t>0%
\mbox{ and }\Vert\mathbf{v}\Vert_{2}=1.
\label{subGau}
\end{equation}
Let $\mathcal{P}_{q}(\tau,c_{n,p})$ denote the set of distributions
of $%
\mathbf{X}_{1}$ satisfying (\ref{subGau}) and with covariance matrix
$\Sigma
\in\mathcal{G}_{q}(c_{n,p})$.

Our technical analysis used in establishing a rate-sharp minimax lower bound
has three major steps. The first step is to reduce the original problem
to a
simpler estimation problem over a carefully chosen subset of the parameter
space without essentially decreasing the level of difficulty. The
second is
to apply the general minimax lower bound technique to this simplified
problem, and the final key step is to bound the total variation affinities
between pairs of mixture distributions with specially designed sparse
covariance matrices. The technical analysis requires ideas that are quite
different from those used in the typical function/sequence estimation
problems.

The minimax upper bound is obtained by studying the risk properties of
thresholding estimators. It will be shown that the optimal rate of
convergence under mean squared spectral norm error is achieved by a
thresholding estimator introduced in \citet{BicLev08N2}. We
write $%
a_{n}\asymp b_{n}$ if there are positive constants $c$ and~$C$ independent
of $n$ such that $c\leq a_{n}/b_{n}\leq C$.
For $1\leq w\leq\infty$, the\vadjust{\goodbreak}
matrix $\ell_{w}$ operator norm of a matrix $A$ is defined by $|\!|\!|
A|\!|\!|_{w}=\max_{\Vert x\Vert_{w}=1}\Vert Ax\Vert_{w}$. The
commonly used spectral norm $|\!|\!|\cdot|\!|\!|$ coincides with the
matrix $\ell_{2}$ operator norm $|\!|\!|\cdot|\!|\!|_{2}$.
(Throughout the paper, we shall write $|\!|\!|\cdot|\!|\!|$ without
a subscript for the matrix spectral norm.) For a symmetric matrix $A$,
it is
known that the spectral norm $|\!|\!|A|\!|\!|$ is equal to the
largest magnitude of the eigenvalues of $A$. Throughout the paper we shall
assume that $1<n^{\beta}\leq p$ for some constants $\beta>1$. Combining
the results given in Sections~\ref{lowbdsec} and~\ref{upperbdsec},
we have
the following optimal rate of convergence for estimating sparse covariance
matrices under the spectral norm.

%
\begin{theorem}
\label{MinimaxOpe} Assume that
%
%
\begin{equation}
c_{n,p}\leq Mn^{{(1-q)}/{2}} ( \log p )^{-{(3-q)}/{2}} \label{condc}
\end{equation}
for $0\le q< 1$. The minimax risk of estimating the covariance matrix
$\Sigma$ under the
spectral norm over the class $\mathcal{P}_{q}(\tau,c_{n,p})$ satisfies
%
%
\begin{equation}
\inf_{\hat{\Sigma}}\sup_{\theta\in\mathcal{P}_{q}(\tau
,c_{n,p})}\mathbb{E}_{\mathbf{X%
}|\theta} |\!|\!|\hat{
\Sigma}-\Sigma|\!|\!|^{2}\asymp c_{n,p}^{2}
\biggl( \frac{\log p}{n} \biggr)^{1-q}+\frac{\log p}{n},
\label{rateOper}
\end{equation}
where $\theta$ denotes a distribution in $\mathcal{P}_{q}(\tau
,c_{n,p})$ with the covariance matrix $\Sigma$.
Furthermore, \eqref{rateOper} holds under the squared $\ell_{w}$
operator norm loss
for all $1\leq w\leq\infty$.
\end{theorem}

We shall focus the discussions on the spectral norm loss. The extension to
the general matrix $\ell_{w}$ operator norm is given in Section \ref
{discussionssec}.
In addition, we also consider
optimal estimation under a class of Bregman matrix divergences which
include Stein's loss, squared Frobenius norm and von Neumann entropy as special
cases. Bregman matrix divergences provide a flexible class of dissimilarity
measures between symmetric matrices and have been used for covariance and
precision matrix estimation as well as matrix approximation problems. See,
for example, \citet{DhiTro07}, Ravikumar et al. (\citeyear
{Ravetal08}) and \citet{KulSusDhi09}. We give a unified result on
the minimax rate of
convergence in Section~\ref{BDsec}.

Besides the sparsity assumption considered in this paper, another commonly
used structural assumption in the literature is that the covariance matrix
is ``bandable'' where the entries decay as they move away from the diagonal.
This is particularly suitable in the setting where the variables
exhibit a
certain ordering structure, which is often the case for time series data.
Various regularization methods have been proposed and studied under this
assumption. \citet{BicLev08N1} proposed a banding estimator and
obtained rate of convergence for the estimator. \citet{CaiZhaZho10}
established the minimax rates of convergence and introduced a rate-optimal
tapering estimator. In particular, \citet{CaiZhaZho10} derived rate
sharp minimax lower bounds for estimating bandable matrices. It should be
noted that the lower bound techniques used there do not lead to a good
result for estimating sparse covariance matrices under the spectral norm.\vadjust{\goodbreak}

The rest of the paper is organized as follows. Section~\ref{Glowbdsec}
introduces a general technical tool for deriving minimax lower bounds
on the
minimax risk. Section~\ref{lowbdsec} establishes the minimax lower bound
for estimating sparse covariance matrices under the spectral norm. The upper
bound is obtained in Section~\ref{upperbdsec} by studying the risk
properties of thresholding estimators.
Section~\ref{BDsec} considers optimal estimation under the Bregman
divergences. A uniform optimal rate of convergence is given for a class of
Bregman divergence losses. Section~\ref{discussionssec} discusses
extensions to estimation under the general $\ell_w$ norm for $1\le w
\le
\infty$ and connections to other related problems including optimal
estimation of sparse precision matrices. The proofs are given in
Section~\ref{proofssec}.

\section{General lower bound for minimax risk}
\label{Glowbdsec}

In this section we develop a new general minimax lower bound technique that
is particularly well suited for treating ``two-directional'' problems such
as estimating sparse covariance matrices. The new method can be viewed
as a
generalization of both Le Cam's method and Assouad's lemma. To help motivate
and understand the new lower bound argument, it is useful to briefly
review Le Cam's method and Assouad's lemma.

Le Cam's method is based on a two-point testing argument and is
particularly well used in estimating linear functionals. See
\citet{LeCam73}
and \citet{DonLiu91}. Let $X$ be an observation from a
distribution $%
\mathbb{P}_{\theta}$ where $\theta$ belongs to a parameter set
$\Theta$.
For two distributions $\mathbb{P}$ and $\mathbb{Q}$ with densities
$p$ and $%
q $ with respect to any common dominating measure $\mu$, the total
variation affinity is given by $\Vert\mathbb{P}\wedge\mathbb
{Q}\Vert=\int
p\wedge q \,d\mu$. Le Cam's method works with a finite parameter set
$\Theta
= \{ \theta_{0},\theta_{1},\ldots,\theta_{D} \} $. Let
$L$ be a
loss function. Define $l_{\min}=\min_{1\leq i\leq D}\inf_{t} [
L (
t,\theta_{0} ) +L ( t,\theta_{i} ) ] $ and
denote $%
\bar{\mathbb{P}}=\frac{1}{D}\sum_{i=1}^{D}\mathbb{P}_{\theta
_{i}}$. Le
Cam's method gives a lower bound for the maximum estimation risk over the
parameter set $\Theta$.

%
\begin{lemma}[(Le Cam)]
\label{LeCam} Let $T$ be any estimator of $\theta$ based on an
observation $%
X$ from a distribution $\mathbb{P}_{\theta}$ with $\theta\in\Theta
= \{ \theta_{0},\theta_{1},\ldots,\theta_{D} \} $, then
%
%
\begin{equation}
\sup_{\theta\in\Theta}\mathbb{E}_{\mathbf{X}|\theta}L ( T,\theta
) \geq
\frac{1}{2}l_{\min}\llVert\mathbb{P}_{\theta
_{0}}\wedge\bar{
\mathbb{P}}\rrVert. \label{LeCamlbd}
\end{equation}
\end{lemma}

Write $\Theta_{1}= \{ \theta_{1},\ldots,\theta_{D} \} $. One
can view the lower bound in (\ref{LeCamlbd}) as obtained from testing the
simple hypothesis $H_{0}\dvtx\theta=\theta_{0}$ against the composite
alternative $H_{1}\dvtx\theta\in\Theta_{1}$.

Assouad's lemma works with a hypercube $\Theta= \{ 0,1 \}^{r}$.
It is based on testing a number of pairs of simple hypotheses and is
connected to multiple comparisons. For a parameter $\theta=(\theta
_{1},\ldots,\theta_{r})$ where $\theta_{i}\in\{0,1\}$, one tests
whether $%
\theta_{i}=0$ or $1$ for each $1\leq i\leq r$ based on the observation $X$.
For each pair of simple hypotheses, there is a certain loss for making an
error in the comparison. The lower bound given by Assouad's lemma is a
combination of losses from testing all pairs of simple hypotheses. Let
%
%
\begin{equation}
H \bigl( \theta,\theta^{\prime} \bigr) =\sum_{i=1}^{r}
\bigl\vert\theta_{i}-\theta_{i}^{\prime}\bigr
\vert\label{H}
\end{equation}
be the Hamming distance on $\Theta$. Assouad's lemma gives a lower bound
for the maximum risk over the hypercube $\Theta$ of estimating an arbitrary
quantity $\psi( \theta) $ belonging to a metric space with
metric $d$.

%
\begin{lemma}[(Assouad)]
\label{Assouad} Let $X\sim\mathbb{P}_{\theta}$ with $\theta\in
\Theta
= \{ 0,1 \}^{r}$, and let $T=T(X)$ be an estimator of $\psi
(\theta
)$ based on $X$. Then for all $s>0$,
%
%
\begin{eqnarray}\label{Assouadlbd}
&&\max_{\theta\in\Theta}2^{s}\mathbb{E}_{\mathbf{X}|\theta
}d^{s}
\bigl( T,\psi( \theta) \bigr)
\nonumber
\\[-8pt]
\\[-8pt]
\nonumber
&&\qquad\geq\min_{H ( \theta
,\theta
^{\prime} ) \geq1}\frac{d^{s} ( \psi( \theta
) ,\psi
( \theta^{\prime} ) ) }{H ( \theta,\theta^{\prime
} ) }\cdot
\frac{r}{2}\cdot\min_{H ( \theta,\theta
^{\prime
} ) =1}\llVert\mathbb{P}_{\theta}\wedge
\mathbb{P}_{\theta
^{\prime}}\rrVert.
\end{eqnarray}
\end{lemma}

We now introduce our new lower bound technique. Again, let $X\sim
\mathbb{P}%
_{\theta}$ where $\theta\in\Theta$. The parameter space $\Theta$ of
interest has a special structure which can be viewed as the Cartesian product
of two components $\Gamma$ and $\Lambda$. For a given positive
integer $r$
and a finite set $B\subset\mathbb{R}^{p}\setminus\{ \mathbf{0}
_{1\times p} \} $, let $\Gamma= \{ 0,1 \}^{r}$ and
$\Lambda
\subseteq B^{r}$. Define
%
%
\begin{equation}
\label{Theta} \Theta=\Gamma\otimes\Lambda= \bigl\{ \theta= (
\gamma,\lambda)
\dvtx\gamma\in\Gamma\mbox{ and }\lambda\in\Lambda\bigr\}.
\end{equation}
In comparison, the standard lower bound arguments work with either
$\Gamma$
or $\Lambda$ alone. For example, Assouad's lemma considers only the
parameter set $\Gamma$ and the Le Cam's method typically applies to a
parameter set like $\Lambda$ with $r=1$. For $\theta=( \gamma
,\lambda)\in
\Theta$, denote the projection of $\theta$ to $\Gamma$ by $\gamma
(\theta) =
\gamma$ and to $\Lambda$ by $\lambda(\theta) = \lambda$.

It is important to understand the structure of the parameter space
$\Theta$%
. One can view an element $\lambda\in\Lambda$ as an $r\times p$ matrix
with each row coming from the set $B$ and view $\Gamma$ as a set of
parameters along the rows indicating whether a given row of $\lambda$ is
present or not. Let $D_{\Lambda}=\operatorname{Card} ( \Lambda
) $.
For a given $a\in\{0,1\}$ and $1\leq i\leq r$, denote $\Theta_{i,a} =
\{\theta\in\Theta\dvtx\gamma_{i}(\theta)=a\}$ where $\theta= (
\gamma
,\lambda) $
and $\gamma_{i}(\theta)$ is the $i$th coordinate of of the first
component of $\theta$. It is easy to see that $\operatorname
{Card}(\Theta_{i,a}) =
2^{r-1} D_\Lambda$. Define the mixture distribution $\bar\mathbb{P}_{i,a}$
by
%
%
\begin{equation}
\bar\mathbb{P}_{i,a}=\frac{1}{2^{r-1}D_{\Lambda}}\sum
_{\theta
\in
\Theta_{i,a}}\mathbb{P}_{\theta} \label{avepi}.
\end{equation}
So $\bar\mathbb{P}_{i,a}$ is the mixture distribution over all
$\mathbb{P}%
_{\theta}$ with $\gamma_{i}(\theta)$ fixed to be $a$ while all other
components of $\theta$ vary over all possible values in $\Theta$.

The following lemma gives a lower bound for the maximum risk over the
parameter set $\Theta$ of estimating a functional $\psi(\theta)$ belonging
to a metric space with metric $d$.

%
\begin{lemma}
\label{AL} For any $s>0$ and any estimator $T$ of $\psi(\theta)$
based on
an observation from the experiment $ \{ \mathbb{P}_{\theta
},\theta\in
\Theta\} $ where $\Theta$ is given in (\ref{Theta}),
%
%
\begin{equation}
\max_{\Theta}2^{s}\mathbb{E}_{\mathbf{X}|\theta}d^{s}
\bigl( T,\psi( \theta) \bigr) \geq\alpha\frac{r}{2}\min_{1\leq
i\leq
r}
\llVert\bar\mathbb{P}_{i,0}\wedge\bar\mathbb{P}_{i,1}
\rrVert, \label{rhsLemma1}
\end{equation}
where $\bar\mathbb{P}_{i,a}$ is defined in equation (\ref{avepi})
and $%
\alpha$ is given by
%
%
\begin{equation}
\alpha=\min_{ \{ (\theta,\theta^{\prime})\dvtx H(\gamma(\theta
),\gamma
(\theta^{\prime}))\geq1 \} }\frac{d^{s}(\psi(\theta),\psi
(\theta^{\prime}))}{H(\gamma(\theta),\gamma(\theta^{\prime}))}.
\label{alpha}
\end{equation}
\end{lemma}

The idea behind this new lower bound argument is similar to the one for
Assouad's lemma, but exists in a more complicated setting. Based on an
observation $%
X\sim\mathbb{P}_{\theta}$ where $\theta=(\gamma,\lambda)\in
\Theta
=\Gamma\otimes\Lambda$, we wish to test whether $\gamma_{i}=0 $ or $1$
for each $1\leq i\leq r$. The first factor $\alpha$ in the lower bound
(\ref%
{rhsLemma1}) is the minimum cost of making an error per comparison. The
second factor $r/2$ is the expected number of errors one makes to
estimate $%
\gamma$ when $\mathbb{P}_{\theta}$ and $\mathbb{P}_{\theta^{\prime}}$
are indistinguishable from each other in the case $H ( \gamma
(\theta
),\gamma(\theta^{\prime}) ) =r$, and the last factor is the lower
bound for the total probability of making type I and type II errors for each
comparison. A major difference is that in this third factor the
distributions $\bar\mathbb{P}_{i,0}$ and $\bar\mathbb{P}_{i,1}
$ are
both complicated mixture distributions instead of the typically simple ones
in Assouad's lemma. This makes the lower bound argument more generally
applicable, while the calculation of the affinity becomes much more
difficult.

In applications of Lemma~\ref{AL}, for a $\gamma=(\gamma_{1},\ldots
,\gamma_{r})\in\Gamma$ where $\gamma_{i}$ takes value $0$ or $1$,
and a
$\lambda
=(\lambda_{1},\ldots,\lambda_{r})\in\Lambda$ where each $\lambda
_{i}\in B$
is a $p$-dimensional nonzero row vector, the element $\theta=(\gamma
,\lambda)\in\Theta$ can be equivalently viewed as an $r\times p$ matrix
%
%
\begin{equation}
\pmatrix{
\gamma_{1}\cdot
\lambda_{1}
\vspace*{2pt}\cr
\gamma_{2}\cdot\lambda_{2}
\vspace*{2pt}\cr
\vdots
\vspace*{2pt}\cr
\gamma_{r}\cdot\lambda_{r}}
, \label{theta}
\end{equation}
where the product $\gamma_{i}\cdot\lambda_{i}$ is taken elementwise:
$%
\gamma_{i}\cdot\lambda_{i}=\lambda_{i}$ if $\gamma_{i}=1$ and the $i$th
row of $\theta$ is the zero vector if $\gamma_{i}=0$. The term $\|
\bar\mathbb{P}_{i,0}\wedge\bar\mathbb{P}_{i,1}\|$ of equation (\ref{rhsLemma1})
is then the lower bound for the total probability of making type I and type
II errors for testing whether or not the $i$th row of $\theta$ is zero.

Note that the lower bound (\ref{rhsLemma1}) reduces to the classical
Assouad lemma when $\Lambda$ contains only one matrix for which every row
is nonzero, and becomes a two-point argument of Le Cam with one point
against a mixture when $r=1$. The proof of this lemma is given in
Section %
\ref{proofssec}. The technical argument is an extension of that of
Assouad's lemma. See \citet{Ass83}, \citet{Yu97} and
\citet{van98}.

The advantage of this method is the ability to break down the lower
bound calculations
for the whole matrix estimation problem into calculations\vadjust{\goodbreak} for individual
rows so that the overall analysis is simplified and more tractable. Although
the tool is introduced here for the purpose of estimating a sparse
covariance matrix, it is of independent interest and is expected to be
useful for solving other matrix estimation problems as well.

Bounding the total variation affinity between two mixture distributions
in (%
\ref{rhsLemma1}) is quite challenging in general. The following well-known
result on the affinity is helpful in some applications. It provides lower
bounds for the affinity between two mixture distributions in terms of the
affinities between simpler distributions in the mixtures.

%
\begin{lemma}
\label{avedis} Let $\overline{\mathbb{P}}_{m}=\sum_{i=1}^{m} w_i
\mathbb{P}%
_{i}$ and $\overline{\mathbb{Q}}_{m}=\sum_{i=1}^{m} w_i \mathbb{Q}_{i}$
where $w_i\ge0$ and $\sum_{i=1}^m w_i =1$. Then
\[
\llVert\overline{\mathbb{P}}_{m}\wedge\overline{\mathbb{Q}}%
_{m}
\rrVert\geq\sum_{i=1}^{m} w_i
\llVert\mathbb{P}_{i}\wedge\mathbb{Q}_{i}\rrVert\geq
\min_{1\leq i\leq m}\llVert\mathbb{P}%
_{i}\wedge
\mathbb{Q}_{i}\rrVert.
\]
\end{lemma}

More specifically, in our construction of the parameter set for establishing
the minimax lower bound, $r$ is the number of possibly nonzero rows in the
upper triangle of the covariance matrix, and $\Lambda$ is the set of
matrices with $r$ rows to determine the upper triangle matrix. Recall that
the projection of $\theta\in\Theta$ to $\Gamma$ is $\gamma(
\theta
) =\gamma= ( \gamma_{i} ( \theta) )_{1\leq
i\leq r}$ and the projection of $\theta$ to $\Lambda$ is $\lambda
(
\theta) =\lambda= ( \lambda_{i} ( \theta)
)_{1\leq i\leq r}$. More generally, for a subset $A\subseteq\{
1,2,\ldots,r \} $, we define a projection of $\theta$ to a
subset of $%
\Gamma$ by $\gamma_{A} ( \theta) = ( \gamma_{i} (
\theta) )_{i\in A}$. A particularly useful example of
set $A$
is
\[
\{ -i \} = \{ 1,\ldots,i-1,i+1,\ldots,r \}
\]
for which $\gamma_{\{-i\}} ( \theta) = ( \gamma_{1} (
\theta) ,\ldots,\gamma_{i-1} ( \theta) ,\gamma_{i+1} ( \theta)
,\gamma_{r} ( \theta)
) $ and in this case for convenience we set $\gamma_{-i}=\gamma_{\{
-i\}}$.
$\lambda_{A} ( \theta) $ and $\lambda_{-i} ( \theta
) $
are defined similarly. We also define the set $\Lambda_{A}= \{
\lambda_{A} ( \theta) \dvtx\theta\in\Theta\} $. A special
case is $%
A=\{-i\}$.

Now we define a subset of $\Theta$ to reduce the problem of estimating
$%
\Theta$ to a problem of estimating $\lambda_{i} ( \theta
) $.
For $a\in\{ 0,1 \} $, $b\in\{ 0,1 \}^{r-1}$
and $%
c\in\Lambda_{-i}\subseteq B^{r-1}$, let
\[
\Theta_{(i,a,b,c)}= \bigl\{ \theta\in\Theta\dvtx\gamma_{i}(
\theta)=a,\gamma_{-i}(\theta)=b\mbox{ and }\lambda_{-i}(
\theta)=c \bigr\}
\]
and $D_{(i,b,c)}=\operatorname{Card}(\Theta_{(i,a,b,c)})$. Note that the
cardinality of $\Theta_{(i,a,b,c)}$ on the right-hand side does not depend
on the value of $a$ due to the Cartesian product structure of $\Theta
=\Gamma
\otimes\Lambda$. Define the mixture distribution
%
%
\begin{equation}
\bar\mathbb{P}_{ ( i,a,b,c ) }=\frac
{1}{D_{(i,b,c)}}\sum
_{\theta
\in\Theta_{(i,a,b,c)}}\mathbb{P}_{\theta}. \label{avepibd}
\end{equation}
In other words, $\bar\mathbb{P}_{(i,a,b,c)}$ is the mixture distribution
over all $\mathbb{P}_{\theta}$ with $\lambda_{i}(\theta)$ varying over
all possible values while all other components of $\theta$ remain
fixed. It
is helpful to observe that when $a=0$, we have $\gamma_{i}(\theta
)\cdot
\lambda_{i}(\theta)=0$ for which $\bar\mathbb{P}_{ (
i,a,b,c )
} $ is degenerate in the sense that it is an average of identical
distributions.

Lemmas~\ref{AL} and~\ref{avedis} together immediately imply the following
result which is based on the total variation affinities between slightly
less complicated mixture distributions. We need to introduce a new
notation $%
\tilde\mathbb{E}_{\theta}$ to denote the average of a function~$g$
over $%
\Theta$, that is,%
\[
\tilde\mathbb{E}_{\theta}g ( \theta) =\sum
_{\theta
\in\Theta}%
\frac{1}{2^{r-1}D_{\Lambda}}g ( \theta) .
\]
The parameter $\theta$ is seen uniformly distributed over $\Theta$. Let
\begin{eqnarray*}
\Theta_{-i} &=& \{ 0,1 \}^{r-1}\otimes\Lambda_{-i}
\\
&=& \bigl\{ ( b,c ) \dvtx\exists\theta\in\Theta\mbox{ such that }
\gamma_{-i}(\theta)=b\mbox{ and }\lambda_{-i}(\theta)=c
\bigr\} ,
\end{eqnarray*}
and an average of $h ( \gamma_{-i},\lambda_{-i} ) $ over
the set $%
\Theta_{-i}$ is defined as follows:%
\[
\tilde\mathbb{E}_{ ( \gamma_{-i},\lambda_{-i} )
}h ( \gamma_{-i},
\lambda_{-i} ) =\sum_{ ( b,c ) \in\Theta
_{-i}}
\frac{%
D_{i,b,c}}{2^{r-1}D_{\Lambda}}h ( b,c ),
\]
where the distribution of $ ( \gamma_{-i},\lambda_{-i} ) $ is
induced by the uniform distribution over $\Theta$.

%
\begin{corollary}
\label{lowbdlemma1} For any $s>0$ and any estimator $T$ of $\psi
(\theta)$
based on an observation from the experiment $ \{ \mathbb
{P}_{\theta
},\theta\in\Theta\} $ where the parameter space $\Theta$ is given
in (\ref{Theta}),
%
%
\begin{eqnarray}
&&\max_{\Theta}2^{s}\mathbb{E}_{\mathbf{X}|\theta}d^{s}
\bigl( T,\psi( \theta) \bigr)
\nonumber
\\
&&\qquad\geq\alpha\frac{r}{2}\min_{i}\tilde\mathbb{E}_{ (
\gamma
_{-i},\lambda_{-i} ) }
\llVert\bar\mathbb{P}_{ (
i,0,\gamma
_{-i},\lambda_{-i} ) }\wedge\bar\mathbb{P}_{ (
i,1,\gamma
_{-i},\lambda_{-i} ) }
\rrVert\label{avelwbd}
\\
&&\qquad\geq\alpha\frac{r}{2}\min_{i}\min_{\gamma_{-i},\lambda
_{-i}}\llVert
\bar\mathbb{P}_{ ( i,0,\gamma_{-i},\lambda_{-i} )
}\wedge\bar\mathbb{P}_{ ( i,1,\gamma_{-i},\lambda_{-i} )
}\rrVert,
\label{minlwbd}
\end{eqnarray}
where $\alpha$ and $\bar\mathbb{P}_{i,a,b,c}$ are defined in
equations (%
\ref{alpha}) and (\ref{avepibd}), respectively.
\end{corollary}

%
\begin{remark}
A key technical step in applying Lemma~\ref{AL} in a typical
application is to show that the affinity $\llVert\bar\mathbb{P}
_{i,0}\wedge\bar\mathbb{P}_{i,1}\rrVert$ is uniformly
bounded away
from $0$ by a constant for all $i$. Then the term $\alpha r$ on the right-hand
side of equation (\ref{rhsLemma1}) in Lemma~\ref{AL} gives the lower
bound for the minimax rate of convergence. As mentioned earlier, the
affinity calculations for two mixture distributions can be very much
involved. Corollary~\ref{lowbdlemma1} gives two lower bounds in terms
of the
affinities. As noted earlier, $\bar\mathbb{P}_{(i,0,\gamma
_{-i},\lambda
_{-i})}$ in the affinity in equations (\ref{avelwbd}) and (\ref
{minlwbd}) is
in fact a single normal distribution, not a mixture. Thus the lower bounds
given in equations~(\ref{avelwbd}) and~(\ref{minlwbd}) require simpler,
although still involved, calculations. In this paper we will apply
equation (%
\ref{avelwbd}), which has an average of affinities on the right-hand
side.
\end{remark}

\section{Lower bound for estimating sparse covariance matrix under the
spectral norm}

\label{lowbdsec}

We now turn to the minimax lower bound for estimating sparse covariance
matrices under the spectral norm. We shall apply the lower bound technique
developed in the previous section to establish rate sharp results. The same
lower bound also holds under the general $\ell_w$ norm for $1\le w \le
\infty$. Upper bounds are discussed in Section~\ref{upperbdsec} and
optimal estimation under Bregman divergence losses is considered in
Section %
\ref{BDsec}.

In this section we shall focus on the Gaussian case and wish to
estimate the
covariance matrix $\Sigma_{p\times p}$ under the spectral norm based
on the
sample $\mathbf{X}_{1},\ldots,\mathbf{X}_{n}\stackrel{\mathrm{i.i.d.}}{\sim
}N(\mu
,\Sigma_{p\times p})$. The parameter space $\mathcal{G}_{q}(c_{n,p})$ for
sparse covariance matrices is defined as in (\ref{sparseparaspace}).
In the
special case of $q=0$, $\mathcal{G}_{0}(c_{n,p})$ contains matrices
with at
most $c_{n,p} + 1$ nonzero elements on each row/column. The parameter
space $%
\mathcal{G}_{q}(c_{n,p})$ also contains the uniformity class $\mathcal
{G}%
_{q}^{\ast}(c_{n,p})$ considered in \citet{BicLev08N2} as a special
case, where $\mathcal{G}_{q}^{\ast}(c_{n,p})$ is defined as, for
$0\leq q<1$,
%
%
\begin{equation}
\mathcal{G}_{q}^{\ast}(c_{n,p})= \biggl\{ \Sigma=(
\sigma_{ij})_{1\leq
i,j\leq p}\dvtx\max_{j\leq p,j\neq i}\sum
_{i\neq j}|\sigma_{ij}|^{q}\leq
c_{n,p} \biggr\} . \label{uniformityclass}
\end{equation}
The columns of $\Sigma\in\mathcal{G}_{q}^{\ast}(c_{n,p})$ are
assumed to
belong to a strong $\ell_{q}$ ball.

We now state and prove the minimax lower bound for estimating a sparse
covariance matrix over the parameter space $\mathcal{G}_{q}(c_{n,p})$ under
the spectral norm. The derivation of the lower bounds relies heavily on the
general lower bound technique developed in the previous section. It also
requires a careful construction of a finite subset of the parameter space
and detailed calculations of an effective lower bound for the total
variation affinities between mixtures of multivariate Gaussian distributions.

%
\begin{theorem}
\label{Operlowerbdthm} Let $\mathbf{X}_{1},\ldots,\mathbf
{X}_{n}\stackrel{%
i.i.d.}{\sim}N(\mu,\Sigma_{p\times p})$. The minimax risk for
estimating the
covariance matrix $\Sigma$ over the parameter space $\mathcal{G}%
_{q}(c_{n,p})$ with $c_{n,p}\leq Mn^{{(1-q)}/{2}} ( \log
p )^{-{(3-q)}/{2}}$ satisfies
%
%
\begin{equation}
\inf_{\hat{\Sigma}}\sup_{\Sigma\in\mathcal
{G}_{q}(c_{n,p})}\mathbb{E}_{\mathbf{X}%
|\Sigma} |\!|\!|\hat{
\Sigma}-\Sigma|\!|\!|^{2}\geq c \biggl( c_{n,p}^{2}
\biggl( \frac{\log p}{n} \biggr)^{1-q}+\frac{\log
p}{n} \biggr)
\label{lowerbd}
\end{equation}
for some constant $c>0$, where $|\!|\!|\cdot|\!|\!|$ denotes the matrix
spectral norm$.$
\end{theorem}

Theorem~\ref{Operlowerbdthm} yields immediately a minimax lower
bound for
the more general subgaussian case under assumption (\ref{subGau}),
\[
\inf_{\hat{\Sigma}}\sup_{\theta\in\mathcal{P}_{q}(\tau
,c_{n,p})}\mathbb{E}_{\mathbf{X%
}|\theta} |\!|\!|\hat{
\Sigma}-\Sigma|\!|\!|^{2}\geq c \biggl( c_{n,p}^{2}
\biggl( \frac{\log p}{n} \biggr)^{1-q}+\frac{\log
p}{n} \biggr) .
\]

It has been shown in \citet{CaiZhaZho10} that%
\[
\inf_{\hat{\Sigma}}\sup_{\theta\in\mathcal{P}_{q}(\tau
,c_{n,p})}\mathbb{E}_{\mathbf{X%
}|\theta} |\!|\!|\hat{
\Sigma}-\Sigma|\!|\!|^{2}\geq c\frac{\log p}{%
n}
\]
by constructing a parameter space with only diagonal matrices. It then
suffices to show that%
\[
\inf_{\hat{\Sigma}}\sup_{\theta\in\mathcal{P}_{q}(\tau
,c_{n,p})}\mathbb{E}_{\mathbf{X%
}|\theta} |\!|\!|\hat{
\Sigma}-\Sigma|\!|\!|^{2}\geq c\cdot c_{n,p}^{2}
\biggl( \frac{\log p}{n} \biggr)^{1-q}
\]
to establish Theorem~\ref{Operlowerbdthm}.

The proof of Theorem~\ref{Operlowerbdthm} contains three major
steps. In
the first step we construct in detail a finite subset $\mathcal{F}_*$
of the
parameter space $\mathcal{G}_q(c_{n,p})$ such that the difficulty of
estimation over $\mathcal{F}_*$ is essentially the same as that of
estimation over $\mathcal{G}_q(c_{n,p})$. The second step is the application
of Lemma~\ref{AL} to the carefully constructed parameter set $\mathcal
{F}_*$. Finally in the third step we calculate the factor $\alpha$
defined in (\ref%
{alpha}) and the total variation affinity between two multivariate normal
mixtures. Bounding the affinity is technically involved. The main ideas of
the proof are outlined here, and detailed proofs of some technical lemmas
used here are deferred to Section~\ref{proofssec}.

\begin{pf*}{Proof of Theorem~\ref{Operlowerbdthm}} The
proof is divided into three main steps.

\textit{Step} 1: \textit{Constructing the parameter set}.
Let $%
r=\lfloor p/2\rfloor$, where $\lfloor x\rfloor$ denotes the largest
integer less than or equal to $x$, and let $B$ be the collection of all row
vectors $b= ( v_{j} )_{1\leq j\leq p}$ such that $v_{j}=0$
for $%
1\leq j\leq p-r$ and $v_{j}=0$ or $1$ for $p-r+1\leq j\leq p$ under the
constraint the total number of 1s is $\llVert b\rrVert_{0}=k$,
where the value of $k$ will be specified later. We shall treat each $%
(b_{1},\ldots,b_{r})\in B^{r}$ as an $r\times p$ matrix with the $i$th row
equal to $b_{i}$.

Set $\Gamma= \{ 0,1 \}^{r}$. Define $\Lambda\subset
B^{r}$ to be
the set of all elements in $B^{r}$ such that each column sum is less
than or
equal to $2k$. For each component $\lambda_{m}$, $1\leq m\leq r$, of $%
\lambda=(\lambda_{1},\ldots,\lambda_{r})\in\Lambda$, define a
$p\times p$
symmetric matrix $A_{m}(\lambda_{m})$ by making the $m$th row of $%
A_{m}(\lambda_{m})$ equal to $\lambda_{m}$, the $m$th column equal to
$%
\lambda_{m}^{T}$ and the rest of the entries $0$. Note that for each $%
\lambda=(\lambda_{1},\ldots,\lambda_{r})\in\Lambda$, each
column/row sum
of the~matrix $\sum_{m=1}^{r}A_{m}(\lambda_{m})$ is less than or
equal to $%
2k$.

Define
%
%
\begin{equation}
\Theta=\Gamma\otimes\Lambda\label{thetaAll},
\end{equation}
and let $\epsilon_{n,p}\in\mathbb{R}$ be fixed. (The exact value of $
\epsilon_{n,p}$ will be chosen later.) For each $\theta=(\gamma
,\lambda
)\in\Theta$ with $\gamma=(\gamma_{1},\ldots,\gamma_{r})\in\Gamma$
and $%
\lambda=(\lambda_{1},\ldots,\lambda_{r})\in\Lambda$, we associate
$\theta$
with a covariance matrix $\Sigma(\theta)$ by
%
%
\begin{equation}
\Sigma(\theta)=I_{p}+\epsilon_{n,p}\sum
_{m=1}^{r}\gamma_{m}A_{m}(
\lambda_{m}). \label{lowbdspace}
\end{equation}
It is easy to see that in the Gaussian case $ |\!|\!|\Sigma_{p\times
p} |\!|\!|\leq\tau$ is a sufficient condition for (\ref{subGau}).
Without loss of generality we assume that $\tau>1$ in the subgaussianity
assumption (\ref{subGau}); otherwise we replace $I_{p}$ in (\ref
{lowbdspace}%
) by $cI_{p}$ with a small constant $c>0$. Finally we define a
collection $%
\mathcal{F}_{\ast}$ of covariance matrices as
%
%
\begin{equation}
\mathcal{F}_{\ast}= \Biggl\{ \Sigma(\theta)\dvtx\Sigma(\theta
)=I_{p}+\epsilon_{n,p}\sum_{m=1}^{r}
\gamma_{m}A_{m}(\lambda_{m}), \theta=(\gamma,
\lambda)\in\Theta\Biggr\} . \label{F*}
\end{equation}
Note that each $\Sigma\in\mathcal{F}_{\ast}$ has value $1$ along
the main
diagonal, and contains an $r\times r$ submatrix, say, $A$, at the upper
right corner, $A^{T}$ at the lower left corner and $0$ elsewhere. Each row
of $A$ is either identically $0$ (if the corresponding $\gamma$ value
is $0$%
) or has exactly $k$ nonzero elements with value $\epsilon_{n,p}$.

We now specify the values of $\epsilon_{n,p}$ and $k$ to ensure
$\mathcal{F}%
_{\ast}\subset\mathcal{G}_{q}(c_{n,p})$. Set $\epsilon
_{n,p}=\upsilon
\sqrt{\frac{\log p}{n}}$ for a fixed small constant $\upsilon$, and
let $%
k=\max( \lceil\frac{1}{2}c_{n,p}\epsilon_{n,p}^{-q} \rceil-1, 0 )$ which
implies
\[
\max_{1\le j\leq p}\sum_{i\neq j}|
\sigma_{ij}|^{q}\leq2k\epsilon_{n,p}^{q}
\leq c_{n,p}.
\]
We require
%
%
\begin{equation}
0<\upsilon< \biggl[ \min\biggl\{ \frac{1}{3},\tau-1 \biggr\}
\frac{1}{M}
\biggr]^{{1}/{(1-q)}} \quad\mbox{and}\quad
\upsilon^{2}<\frac{\beta-1%
}{54\beta}. \label{v}
\end{equation}
Note that $\epsilon_{n,p}$ and $k$ satisfy
%
%
\begin{equation}
2k\epsilon_{n,p}\leq c_{n,p}\epsilon_{n,p}^{1-q}
\leq M\upsilon^{1-q}<\min\bigl\{ \tfrac{1}{3},\tau-1 \bigr\}
\label{1normbound}
\end{equation}
and consequently every $\Sigma(\theta)$ is diagonally dominant and
positive definite, and $|\!|\!|\Sigma(\theta)|\!|\!|\leq|\!|\!|
\Sigma
(\theta)|\!|\!|_{1}\leq2k\epsilon_{n,p}+1<\tau$. Thus we have
$\mathcal{F}%
_{\ast}\subset\mathcal{G}_{q}(c_{n,p})$, and the subgaussianity
assumption (%
\ref{subGau}) is satisfied.

\textit{Step} 2: \textit{Applying the general lower bound
argument}. Let $\mathbf{X}_{1},\ldots,\mathbf{X}_{n}\stackrel
{\mathrm{i.i.d.}}{\sim}%
N ( 0,\Sigma(\theta) ) $ with $\theta\in\Theta$ and
denote the
joint distribution by $\mathbb{P}_{\theta}$. Applying Lemma~\ref{AL}
to the
parameter space $\Theta$ with $s=2$, we have
%
%
\begin{equation}
\inf_{\hat{\Sigma}}\max_{\theta\in\Theta}2^{2}\mathbb
{E}_{\mathbf{X}%
|\theta} \bigl|\!\bigl|\!\bigl|\hat{\Sigma}-\Sigma(\theta) \bigr|\!\bigr|\!
\bigr|^{2}\geq
\alpha\cdot\frac{r}{2}\cdot\min_{1\leq i\leq r}\llVert\bar\mathbb{P}%
_{i,0}\wedge\bar\mathbb{P}_{i,1}\rrVert,
\label{lowerbound*}
\end{equation}
where
%
%
\begin{equation}
\alpha\equiv\min_{ \{ (\theta,\theta^{\prime})\dvtx H(\gamma
(\theta
),\gamma(\theta^{\prime}))\geq1 \} }\frac{ |\!|\!|
\Sigma(\theta
)-\Sigma(\theta^{\prime}) |\!|\!|^{2}}{H(\gamma(\theta
),\gamma
(\theta^{\prime}))} \label{alpha1},
\end{equation}
and $\bar\mathbb{P}_{i,0}$ and $\bar\mathbb{P}_{i,1}$ are
defined as in (%
\ref{avepi}).

\textit{Step} 3: \textit{Bounding the affinity and per
comparison loss}. We shall now bound the two factors $\alpha$ and $%
\min_{i}\llVert\bar\mathbb{P}_{i,0}\wedge\bar\mathbb{P}%
_{i,1}\rrVert$ in (\ref{lowerbound*}). This is done separately
in the
next two lemmas which are proved in detail in Section~\ref{proofssec}.
Lemma~\ref{dffbd} gives a lower bound to the per comparison loss, and
it is
easy to prove.

%
\begin{lemma}
\label{dffbd} For $\alpha$ defined in equation (\ref{alpha1}) we have
\[
\alpha\geq\frac{(k\epsilon_{n,p})^{2}}{p}.
\]
\end{lemma}

The key technical difficulty is in bounding the affinity between the
Gaussian mixtures. The proof is quite involved.

%
\begin{lemma}
\label{affbd} Let $\mathbf{X}_{1},\ldots,\mathbf{X}_{n}\stackrel
{i.i.d.}{\sim}%
N ( 0,\Sigma( \theta) ) $ with $\theta\in
\Theta$
defined in equation~(\ref{thetaAll}), and denote the joint
distribution by $%
\mathbb{P}_{\theta}$. For $a\in\{0,1\}$ and $1\leq i\leq r$, define $
\bar\mathbb{P}_{i,a}$ as in (\ref{avepi}). Then there exists a
constant $%
c_{1}>0$ such that
\[
\min_{1\leq i\leq r}\llVert\bar\mathbb{P}_{i,0}\wedge\bar\mathbb{P}%
_{i,1}\rrVert\geq c_{1}.
\]
\end{lemma}

Finally, the minimax lower bound for estimation over $\mathcal{G}%
_{q}(c_{n,p})$ is obtained by putting together the bounds given in
Lemmas %
\ref{dffbd} and~\ref{affbd},
\begin{eqnarray*}
\inf_{\hat{\Sigma}}\sup_{\Sigma\in\mathcal
{G}_{q}(c_{n,p})}\mathbb{E}_{\mathbf{X}%
|\Sigma} |\!|\!|\hat{
\Sigma}-\Sigma|\!|\!|^{2} &\geq& \inf_{\hat{\Sigma}}
\max_{\Sigma
(\theta)\in\mathcal{F}_{\ast}}\mathbb{E}_{\mathbf{X}|\theta
} \bigl|\!\bigl|\!\bigl|\hat{\Sigma}-\Sigma(
\theta) \bigr|\!\bigr|\!\bigr|^{2}\\
&\geq&\frac{ (
k\epsilon_{n,p} )^{2}}{p}\cdot\frac{r}{8}\cdot
c_{1}
\\
&\geq&c_{2}c_{n,p}^{2} \biggl( \frac{\log p}{n}
\biggr)^{1-q}
\end{eqnarray*}
for some constant $c_{2}>0$.
\end{pf*}
%
\begin{remark}
It is easy to check that the proof of Theorem \ref
{Operlowerbdthm} also yields a lower
bound for estimation under the general matrix $\ell_{w}$ operator norm
for any $1\leq w\leq\infty,$
\[
\inf_{\hat{\Sigma}}\sup_{\theta\in\mathcal{P}_{q}(\tau
,c_{n,p})}\mathbb{E}_{\mathbf{X%
}|\theta} |\!|\!|\hat{
\Sigma}-\Sigma|\!|\!|_{w}^{2}\geq c \biggl(
c_{n,p}^{2} \biggl( \frac{\log p}{n} \biggr)^{1-q}+
\frac{\log
p}{n} \biggr)
\]
by applying Lemma~\ref{AL} with $s=1$.
\end{remark}

\section{Minimax upper bound under the spectral norm}
\label{upperbdsec}

Section~\ref{lowbdsec} developed a minimax lower bound for estimating a
sparse covariance matrix under the spectral norm over $\mathcal{G}%
_{q}(c_{n,p})$. In this section we shall show that the lower bound is
rate-sharp and therefore establish the optimal rate of convergence. To
derive a minimax upper bound, we shall consider the properties of a
thresholding estimator introduced in \citet{BicLev08N2}. Given a
random sample $\{\mathbf{X}_{1},\ldots,\mathbf{X}_{n}\}$ of $p$-variate
observations drawn from a distribution in $\mathcal{P}_{q}(\tau,c_{n,p})$,
the sample covariance matrix is
\[
\frac{1}{n-1}\sum_{l=1}^{n} (
\mathbf{X}_{l}-\mathbf{\bar{X}} ) ( \mathbf{X}_{l}-
\mathbf{\bar{X}} )^{T},
\]
which is an unbiased estimate of $\Sigma$, and the maximum likelihood
estimator of $\Sigma$ is
%
%
\begin{equation}
\Sigma^{\ast}=\bigl(\sigma_{ij}^{\ast}
\bigr)_{1\leq i,j\leq p}=\frac{1}{n}%
\sum
_{l=1}^{n} ( \mathbf{X}_{l}-\mathbf{
\bar{X}} ) ( \mathbf{X}%
_{l}-\mathbf{\bar{X}}
)^{T} \label{MLE}
\end{equation}
when $\mathbf{X}_{l}$'s are normally distributed. These two estimators are
close to each other for large $n$. We shall construct estimators of the
covariance matrix $\Sigma$ by thresholding the maximum likelihood estimator
$\Sigma^{\ast}$.

Note that the subgaussianity condition (\ref{subGau}) implies%
\[
|\!|\!|\Sigma|\!|\!|=\sup_{\mathbf{v}\dvtx\llVert
\mathbf{v}%
\rrVert=1} \operatorname{ Var} \bigl[
\mathbf{v}^{T}(\mathbf{X}_{1}-\mathbb{E}
\mathbf{X}_{1}) \bigr] \leq\int_{0}^{\infty
}e^{-x/ ( 2\tau
) }\,dx=2
\tau.
\]
Then the empirical covariance $\sigma_{i,j}^{\ast}$ satisfies the
following large deviation result that there exist constants $C_{1}>0$
and $%
\gamma>0$ such that
%
%
\begin{equation}
\mathbb{P} \bigl( \bigl\vert\sigma_{ij}^{\ast}-
\sigma_{ij}\bigr\vert>t \bigr) \leq C_{1}\exp\biggl( -{
\frac{8}{\gamma^{2}}}nt^{2} \biggr) \label{emsigmaTail}
\end{equation}
for $\vert t\vert\leq\delta$, where $C_{1},$ $\gamma$
and $%
\delta$ are constants and depend only on $\tau$. See Saulis and
Statulevi\v{c}ius (\citeyear{SauSta91}) and \citet{BicLev08N1}.
Inequality (\ref%
{emsigmaTail}) implies $\sigma_{ij}^{\ast}$ behaves like a subgaussian
random variable. In particular for $t=\gamma\sqrt{\frac{\log p}{n}}$ we
have
%
%
\begin{equation}
\mathbb{P} \bigl( \bigl\vert\sigma_{ij}^{\ast}-
\sigma_{ij}\bigr\vert>t \bigr) \leq C_{1}p^{-8}
.\label{sigmalijTail}
\end{equation}
Define the thresholding estimator $\hat{\Sigma}= ( \hat{\sigma}
_{ij} )_{p\times p}$ by
%
%
\begin{equation}
\hat{\sigma}_{ij}=\sigma_{ij}^{\ast}\cdot I \biggl(
\bigl|\sigma_{ij}^{\ast
}\bigr|\geq\gamma\sqrt{\frac{\log p}{n}} \biggr)
. \label{sigmahatij}
\end{equation}
This thresholding estimator was first proposed in \citet{BicLev08N2}
in which a rate of convergence of the loss function in probability was given
over the uniformity class $\mathcal{G}_{q}^{\ast}(c_{n,p})$. Here we
provide an upper bound for mean squared spectral norm error over the
parameter space $\mathcal{G}_{q}(c_{n,p})$.

Throughout the rest of the paper we denote by $C$ a generic positive
constant which may vary from place to place. The following theorem shows
that the thresholding estimator defined in (\ref{sigmahatij}) is rate
optimal over the parameter space {$\mathcal{G}_{q}(c_{n,p})$}.

%
\begin{theorem}
\label{MinimaxOpe1} The thresholding estimator $\hat{\Sigma}$ given
in (\ref%
{sigmahatij}) satisfies, for some constant $C>0$,
%
%
\begin{equation}
\sup_{\theta\in\mathcal{P}_{q}(\tau,c_{n,p})}\mathbb{E}_{\mathbf
{X}|\theta
} |\!|\!|\hat{\Sigma}-\Sigma
|\!|\!|^{2}\leq C \biggl[ c_{n,p}^{2} \biggl(
\frac{\log p}{n} \biggr)^{1-q}+\frac{\log
p}{n} \biggr] .
\label{rateOper2}
\end{equation}
Consequently, the minimax risk of estimating the sparse covariance
matrix $%
\Sigma$ over {$\mathcal{G}_{q}(c_{n,p})$} satisfies
%
%
\begin{equation}
\inf_{\hat{\Sigma}}\sup_{\theta\in\mathcal{P}_{q}(\tau
,c_{n,p})}\mathbb{E}_{\mathbf{X%
}|\theta} |\!|\!|\hat{
\Sigma}-\Sigma|\!|\!|^{2}\asymp c_{n,p}^{2}
\biggl( \frac{\log p}{n} \biggr)^{1-q}+\frac{\log
p}{n}.
\label{gminimax}
\end{equation}
\end{theorem}

%
\begin{remark}
A similar argument to the proof of equation (\ref{rateOper2}) in
Section~\ref{uppersec} yields the following upper bound for estimation under the
matrix $\ell_1$ norm:
\[
\sup_{\theta\in\mathcal{P}_{q}(\tau,c_{n,p})}\mathbb{E}_{\mathbf
{X}|\theta
} |\!|\!|\hat{\Sigma}-\Sigma
|\!|\!|_{1}^{2}\leq C \biggl[ c_{n,p}^{2}
\biggl( \frac{\log p}{n} \biggr)^{1-q}+\frac{\log
p}{n} \biggr] .
\]
\end{remark}

Theorem~\ref{MinimaxOpe1} shows that the optimal rate of convergence for
estimating a sparse covariance matrix over $\mathcal{G}_{q}(c_{n,p})$ under
the squared spectral norm is $c_{n,p}^{2} ( \frac{\log
p}{n} )^{1-q}$. In \citet{BicLev08N2} the uniformity class
$\mathcal{G}
_{q}^{\ast}(c_{n,p})$ defined in (\ref{uniformityclass}) was considered.
We shall now show that the same minimax rate of convergence holds for
estimation over $\mathcal{G}_{q}^{\ast}(c_{n,p})$. It is easy to
check in
the proof of the lower bound that for every $\Sigma\in\mathcal
{F}_{\ast}$
defined in (\ref{F*}), we have
\[
\max_{1\leq j\leq p}\sum_{i\neq j}\vert
\sigma_{ij}\vert^{q}\leq2k\epsilon_{n,p}^{q}
\le c_{n,p}
\]
and consequently $\mathcal{F}_{\ast}\subset\mathcal{G}_{q}^{\ast
}(c_{n,p}) $. Thus the lower bound established for $\mathcal{F}_{\ast}$
automatically yields a lower bound for $\mathcal{G}_{q}^{\ast}(c_{n,p})$.
On the other hand, since a strong $\ell_{q}$ ball is always contained
in a
weak $\ell_{q}$ ball by the Markov inequality, the upper bound in
equation (%
\ref{rateOper2}) for the parameter space ${\mathcal{G}_{q}}$ also
holds for $%
\mathcal{G}_{q}^{\ast}(c_{n,p})$. Let $\mathcal{P}_{q}^{\ast}(\tau
,c_{n,p})$ denote the set of distributions of $\mathbf{X}_{1}$
satisfying (%
\ref{subGau}) and with covariance matrix $\Sigma\in\mathcal
{G}_{q}^{\ast
}(c_{n,p})$. Then we have the following result.

%
\begin{proposition}
\label{BLprop} The minimax risk for estimating the covariance matrix under
the spectral norm over the uniformity class $\mathcal{G}_{q}^{\ast
}(c_{n,p}) $ satisfies
\[
\inf_{\hat{\Sigma}}\sup_{\theta\in\mathcal{P}_{q}^{\ast}(\tau
,c_{n,p})}\mathbb{E}_{%
\mathbf{X}|\theta} |\!|\!|\hat{
\Sigma}-\Sigma|\!|\!|^{2}\asymp c_{n,p}^{2}
\biggl( \frac{\log p}{n} \biggr)^{1-q}+\frac{\log p}{n}.
\]
\end{proposition}

The thresholding estimator $\hat\Sigma$ defined by (\ref
{sigmahatij}) is
positive definite with high probability, but it is not guaranteed to be
positive definite. A simple additional step can make the final estimator
positive semi-definite and achieve the optimal rate of convergence. Write
the eigen-decomposition of $\hat{\Sigma}$ as
\[
\hat{\Sigma}=\sum_{i=1}^{p}\hat{
\lambda}_{i}v_{i}v_{i}^{T},
\]
where $\hat{\lambda}_{i}$'s and $v_i$'s are the eigenvalues and eigenvectors
of $\hat{\Sigma}$, respectively. Let $\hat{\lambda}_{i}^{+}=\max
(\hat
\lambda_i, 0)$ be the positive part of $\hat{\lambda}_{i}$ and define
\[
\hat{\Sigma}^{+}=\sum_{i=1}^{p}
\hat{\lambda}_{i}^{+}v_{i} v_{i}^{T}.
\]
Then%
\begin{eqnarray*}
\bigl|\!\bigl|\!\bigl|\hat{\Sigma}^{+}-\Sigma\bigr|\!\bigr|\!
\bigr|&\leq&
\bigl|\!\bigl|\!\bigl|\hat{
\Sigma}%
^{+}-\hat{\Sigma} \bigr|\!\bigr|\!\bigr|+ |\!|\!|\hat{\Sigma
}-\Sigma
|\!|\!|\leq\max_{i:\hat{\lambda}_{i}\leq0}\vert\hat{\lambda
}_{i}\vert+
|\!|\!|\hat{\Sigma}-\Sigma|\!|\!|
\\
&\leq&\max_{i:\hat{\lambda}_{i}\leq0}\vert\hat{\lambda}_{i}-
\lambda_{i}\vert+ |\!|\!|\hat{\Sigma}-\Sigma|\!|\!|\leq2 |\!|\!|
\hat{\Sigma}-\Sigma|\!|\!|.
\end{eqnarray*}
%
The resulting estimator $\hat{\Sigma}^{+}$ is positive semi-definite and
attains the same rate as the original thresholding estimator $\hat
\Sigma$.
This method can be applied to the tapering estimator in \citet
{CaiZhaZho10} as well to make the estimator positive semi-definite,
while still
achieving the optimal rate.

\section{Optimal estimation under Bregman divergences}
\label{BDsec}

We have so far focused on the optimal rate of convergence under the spectral
norm. In this section we turn to minimax estimation of sparse covariance
matrices under a class of Bregman divergence losses which include
Stein's loss, Frobenius norm and von Neumann's entropy as special cases.
Bregman matrix divergences have been used for matrix estimation and matrix
approximation problems; see, for example, \citet{DhiTro07},
Ravikumar et
al. (\citeyear{Ravetal08}) and \citet{KulSusDhi09}. In this
section we
establish the optimal rate of convergence uniformly for a class of Bregman
divergence losses.

Bregman (\citeyear{Bre67}) introduced the Bregman divergence as a
dissimilarity measure
between vectors,
\[
D_{\phi} ( \mathbf{x},\mathbf{y} ) =\phi( \mathbf{x} ) -\phi(
\mathbf{y} ) - \bigl( \nabla\phi( \mathbf{y} ) \bigr)^{T} (
\mathbf{x}-\mathbf{y} ),
\]
where $\phi$ is a differentiable, real-valued, and strictly convex function
defined over a convex set in a Euclidean space $\mathbb{R}^{m}$, and $%
\nabla\phi$ is the gradient of $\phi$. The well-known Mahalanobis distance
is a Bregman divergence. This concept can be naturally extended to the space
of real and symmetric matrices as
\[
D_{\phi} ( X,Y ) =\phi( X ) -\phi( Y ) -\operatorname{tr} \bigl[
\bigl( \nabla
\phi(
Y ) \bigr)^{T} ( X-Y ) %
\bigr],
\]
where $X$ and $Y$ are real symmetric matrices, and $\phi$ is a
differentiable strictly convex function over the space. See \citet
{CenZen97} and \citet{KulSusDhi09}. A particularly
interesting class of $\phi$ is
%
%
\begin{equation}
\phi( X ) =\sum_{i=1}^{p}\varphi(
\lambda_{i} ), \label{SpectralBD}
\end{equation}
where $\lambda_{i}$'s are the eigenvalues of $X$, and $\varphi$ is a
differentiable, real-valued, and strictly convex function over a convex set
in $\mathbb{R}$.\vadjust{\goodbreak} See \citet{DhiTro07} and \citet
{KulSusDhi09}. Examples of this class of Bregman divergences include:

\begin{itemize}
\item$\varphi( \lambda) =$ $-\log\lambda,$ or
equivalently $%
\phi( X ) =-\log\det( X ) $. The corresponding
Bregman divergence can be written as%
\[
D_{\phi} ( X,Y ) =\operatorname{tr} \bigl( XY^{-1} \bigr) -\log
\det\bigl(
XY^{-1} \bigr) -p,
\]
which is often called Stein's loss in the statistical literature.

\item$\varphi( \lambda) =$ $\lambda\log\lambda
-\lambda,$
or equivalently $\phi( X ) =\operatorname{tr} ( X\log X-X ) $,
where $X$
is positive definite such that $\log X$ is well defined. The corresponding
Bregman divergence is the von Neumann divergence%
\[
D_{\phi} ( X,Y ) =\operatorname{tr} ( X\log X-X\log Y-X+Y ).
\]

\item$\varphi( \lambda) =$ $\lambda^{2},$ or
equivalently $%
\phi( X ) =\operatorname{tr} ( X^{2} ) $. The resulting Bregman
divergence is the squared Frobenius norm%
\[
D_{\phi} ( X,Y ) =\operatorname{tr} \bigl[ ( X-Y )^{2} \bigr] = |\!
|\!|X-Y
|\!|\!|_{F}^{2}=\sum_{i,j} (
x_{ij}-y_{ij} )^{2}
\]
for $X= ( x_{ij} )_{1\leq i,j\leq p}$ and $Y= (
y_{ij} )_{1\leq i,j\leq p}$.
\end{itemize}

Define a class $\Psi$ of functions $\varphi$ satisfying
the following conditions:

\begin{longlist}[(1)]
\item[(1)] $\varphi$ is twice differentiable, real-valued and
strictly convex
over $\lambda\in( 0,\infty) $;

\item[(2)] $\vert\varphi( \lambda) \vert
\leq C\lambda^{r}$ for some $C>0$ and some real number $r$ uniformly
over $\lambda
\in
( 0,\infty) $;

\item[(3)] For every positive constants $\epsilon_{2}$ and $M_{2}$
there are
some positive constants $c_{L}$ and $c_{u}$ depending on $\epsilon
_{2}$ and
$M_{2}$ such that $c_{L}\leq\varphi^{^{\prime\prime}} (
\lambda
) \leq c_{u}$ for all $\lambda\in[ \epsilon_{2},M_{2} ] $.
\end{longlist}

In this paper, we shall consider the following class of Bregman divergences:
%
%
\begin{equation}
\Phi= \Biggl\{ \phi( \Sigma) =\sum_{i=1}^{p}
\varphi( \lambda_{i} ) \dvtx\varphi\in\Psi\Biggr\}.
\label{assumpphi}
\end{equation}
It is easy to see that Stein's loss, von Neumann's divergence and the squared
Frobenius norm are in this class.

Let $\epsilon_{1}>0$ be a positive constant. Let $\mathcal
{P}_{q}^{B}(\tau
,c_{n,p})$ denote the set of distributions of $\mathbf{X}_{1}$
satisfying (%
\ref{subGau}) and with covariance matrix
\[
\Sigma\in\mathcal{G}_{q}^{B}(c_{n,p})=
\mathcal{G}_{q}(c_{n,p})\cap\{ \Sigma\dvtx
\lambda_{\mathrm{min}}\geq\epsilon_{1} \} .
\]
Here $\lambda_{\mathrm{min}}$ denotes the minimum eigenvalue of
$\Sigma$.
The assumption that all eigenvalues are bounded away from $0$ is necessary
when $\varphi( \lambda) $ is not well defined at~$0$. An
example is the Stein loss where $\varphi( \lambda)
=-\log
\lambda$. Under this assumption all losses $D_{\phi}$ are equivalent to
the squared Frobenious norm.

The following theorem gives a unified result on the minimax rate of
convergence for estimating the covariance matrix over the parameter
space $%
\mathcal{P}_{q}^{B}(\tau,c_{n,p})$ for all Bregman divergences $\phi
\in
\Phi$ defined in (\ref{assumpphi}).

%
\begin{theorem}
\label{MinimaxBD}
Assume that $c_{n,p}\leq Mn^{{(1-q)}/{2}}(\log p)^{-{{(3-q)}/{2}}}$
for some $M>0$ and $0\le q< 1$. The minimax risk over
$\mathcal{P}%
_{q}^{B}(\tau,c_{n,p})$ under the loss function
\[
\mathrm{L}_{\phi} ( \hat{\Sigma},\Sigma) =\frac
{1}{p}D_{\phi
}
( \hat{\Sigma},\Sigma)
\]
for all Bregman divergences $\phi\in\Phi$ defined in (\ref{assumpphi})
satisfies
%
%
\begin{equation}
\inf_{\hat{\Sigma}} \sup_{\phi\in\Phi}\sup_{\theta\in\mathcal
{P}_{q}^{B}(\tau
,c_{n,p})}
\mathbb{E}_{\mathbf{X}|\theta}\mathrm{L}_{\phi} ( \hat{\Sigma%
},
\Sigma) \asymp c_{n,p} \biggl( \frac{\log p}{n} \biggr)^{1-{{q}/{2%
}}}+
\frac{1}{n}. \label{BDRisk}
\end{equation}
\end{theorem}

Note that Theorem~\ref{MinimaxBD} gives the minimax rate of convergence
uniformly under all Bregman divergences defined in (\ref{assumpphi}).
For an
individual Bregman divergence loss, the condition that all eigenvalues are
bounded away from $0$ is not needed if the function $\varphi$ is well
behaved at $0$. For example, such is the case for the Frobenius norm.

The optimal\vspace*{1pt} rate of convergence is attained by a modified thresholding
estimator. Let $\hat{\Sigma}=(\hat{\sigma}_{ij})_{1\leq i,j\leq p}$
be the
thresholding estimator given in (\ref{sigmahatij}). Define the final
estimator of $\Sigma$ by
%
%
\begin{equation}\qquad
\hat{\Sigma}_{B}=\cases{
\hat{
\Sigma}, &\quad$\mbox{if } \displaystyle\frac{1}{\max\{ \log
n,\log
p \} }%
\leq\lambda_{\min}(
\hat{\Sigma}) \leq\max\{ \log n,\log p \},$
\vspace*{2pt}\cr
I, & \quad$\mbox{otherwise.}$}
\label{sigmahatStein}
\end{equation}
It will be proved in Section~\ref{Bregmanproofsec} that the
estimator $%
\hat{\Sigma}_{B}$ given in (\ref{sigmahatStein}) is rate optimal uniformly\vspace*{1pt}
under all Bregman divergences satisfying (\ref{assumpphi}). Note that the
modification of $\hat{\Sigma}$ given in (\ref{sigmahatStein}) is needed.
Without it, the loss $\mathrm{L}_{\phi}(\hat{\Sigma},\Sigma)$ may
not be
well behaved under some Bregman divergences such as Stein's loss and von
Neumann's divergence.

%
\begin{remark}
Let $\mathcal{P}_{q}^{\ast B}(\tau,c_{n,p})$ denote the set of
distributions of $\mathbf{X}_{1}$ satisfying~(\ref{subGau}) and with
covariance matrix $\Sigma\in\mathcal{G}_{q}^{\ast
B}(c_{n,p})=\mathcal{G}%
_{q}^{\ast}(c_{n,p})\cap\{ \Sigma\dvtx\lambda_{min}\geq\epsilon
_{1} \} $. Then under the same conditions as in Theorem~\ref{MinimaxBD},
\[
\inf_{\hat{\Sigma}}\sup_{\phi\in\Phi}\sup_{\theta\in\mathcal
{P}_{q}^{\ast B}(\tau
,c_{n,p})}\mathbb{E}_{\mathbf{X}|\theta}
\mathrm{L}_{\phi} ( \hat{\Sigma%
},\Sigma) \asymp
c_{n,p} \biggl( \frac{\log p}{n} \biggr)^{1-{{q}/{2%
}}}+
\frac{1}{n}.
\]
\end{remark}

\section{Discussions}
\label{discussionssec}

The focus of this paper is mainly on the optimal estimation under the
spectral norm. However, both the lower and upper bounds can be easily
extended to the general matrix $\ell_w$ norm for $1\le w \le\infty$ by
using similar arguments given in Sections~\ref{lowbdsec} and~\ref{upperbdsec}.\vadjust{\goodbreak}

%
\begin{theorem}
\label{MinimaxOper} Under the assumptions in Theorem \ref{MinimaxOpe}, the
minimax risk of estimating the covariance matrix $\Sigma$ under the
matrix $%
\ell_{w}$-norm for $1\leq w\leq\infty$ over the class $\mathcal{P}%
_{q}(\tau,c_{n,p})$ satisfies
%
%
\begin{equation}
\inf_{\hat{\Sigma}}\sup_{\theta\in\mathcal{P}_{q}(\tau
,c_{n,p})}\mathbb{E}_{\mathbf{X%
}|\theta} |\!|\!|\hat{
\Sigma}-\Sigma|\!|\!|_w^{2}\asymp c_{n,p}^{2}
\biggl( \frac{\log p}{n} \biggr)^{1-q}+\frac
{\log p}{n%
}.
\label{lw-opt}
\end{equation}
Moreover, the thresholding estimator $\hat{\Sigma}$ defined in (\ref%
{sigmahatij}) is rate-optimal.
\end{theorem}

As noted in Section~\ref{lowbdsec}, a rate-sharp lower bound for the
minimax risk under the $\ell_{w}$ norm can be obtained by using essentially
the same argument with the same parameter space $\mathcal{F}_{\ast}$
and a
slightly modified version of Lemma~\ref{dffbd}. The upper bound can be
proved by applying the Riesz--Thorin interpolation theorem, which yields
$ |\!|\!|A |\!|\!|_{w}\leq\max\{ { |\!|\!|A
|\!|\!|_{1},} |\!|\!|A |\!|\!|_{2},|\!|\!|
A|\!|\!|_{\infty} \} ${ for all }$w\in[ 1,\infty
) $, and by using the facts $|\!|\!|A|\!|\!|_{1}=|\!|\!|A|\!|\!
|_{\infty}$ and $ |\!|\!|A |\!|\!|_{2}\leq|\!|\!|
A |\!|\!|_{1}$, when $A$ is symmetric. In Section \ref
{upperbdsec} we
have in fact established the same risk bound for both the spectral norm and
matrix $\ell_{1}$-norm.

The spectral norm of a matrix depends on the entries in a subtle way
and the
``interactions'' among different rows/columns must be taken into account.
The lower bound argument developed in this paper is aimed at treating
``two-directional'' problems by mixing over both rows and columns. It
can be
viewed as a simultaneous application of Le Cam's method in one
direction and
Assouad's lemma in another. In contrast, for sequence estimation
problems, we
typically need one or the other, but not both at the same time. The lower
bound techniques developed in this paper can be used to solve other matrix
estimation problems. For example, \citet{CaiLiuZho} applied the
general lower bound argument to the problem of estimating sparse precision
matrices under the spectral norm and established the optimal rate of
convergence. This problem is closely connected to graphical model selection.
The derivations of both the lower and upper bounds are involved. For reasons
of space, we shall report the results elsewhere.

In this paper we also developed a unified result on the minimax rate of
convergence for estimating sparse covariance matrices under a class of
Bregman divergence losses which include the commonly used Frobenius
norm as
a special case. The optimal rate of convergence given in Theorem \ref%
{MinimaxBD} is identical to the minimax rate for estimating a
row/column as
a vector with the weak $\ell_{q}$ ball constraint under the squared error
loss. Our result shows that this class of Bregman divergence losses are
essentially the same and thus can be studied simultaneously in terms of the
minimax rate of convergence.

Estimating a sparse covariance matrix is intrinsically a heteroscedastic
problem in the sense that the variances of the entries of the sample
covariance matrix are not equal and can vary over a wide range. A natural
approach is to adaptively threshold the entries according to their
individual variabilities. \citet{CaiLiu11} considered such an adaptive\vadjust{\goodbreak}
approach for estimation over the weighted $\ell_q$ balls which
contains the
strong $\ell_q$ balls as subsets. The lower bound given in Proposition
\ref{BLprop} in the present paper immediately yields a lower bound for
estimation over the weighted $\ell_q$ balls. A data-driven thresholding
procedure was introduced and shown to adaptively achieve the optimal
rate of
convergence over a large collection of the weighted $\ell_q$ balls
under the
spectral norm. In contrast, universal thresholding estimators are
sub-optimal over the same parameter spaces.

In addition to the hard thresholding estimator used in \citet
{BicLev08N2}, \citet{RotLevZhu09} considered a class of thresholding
rules with more general thresholding functions, including soft thresholding
and adaptive Lasso. It is straightforward to show that these thresholding
estimators with the same choice of threshold level used in (\ref{sigmahatij}%
) also attains the optimal rate of convergence over the parameter space
$%
\mathcal{G}_{q}(c_{n,p})$ under mean squared spectral norm error as
well as
under the class of Bregman divergence losses considered in
Section~\ref{BDsec} with the same modification as in (\ref{sigmahatStein}). Therefore,
the choice of the thresholding function is not important as far as the rate
optimality is concerned.

\section{Proofs}
\label{proofssec}

In this section we prove the general lower bound result given in Lemma
\ref%
{AL}, Theorems~\ref{MinimaxOpe1} and~\ref{MinimaxBD} as well as some
of the
important technical lemmas used in the proof of Theorem~\ref{Operlowerbdthm} given in Section~\ref{lowbdsec}. The proofs of a few
technical results used in this section are deferred to the supplementary material [\citet{supp}].
Throughout this section, we denote by $C$ a generic constant that may vary
from place to place.

\subsection{\texorpdfstring{Proof of Lemma \protect\ref{AL}}{Proof of Lemma 3}}

We first bound the maximum risk by the average over the whole parameter set,
%
%
\begin{eqnarray}\label{MaxAve}
\max_{\Theta}2^{s}\mathbb{E}_{\mathbf{X}|\theta}d^{s}
\bigl( T,\psi( \theta) \bigr) &\geq&\frac{1}{2^{r}D_{\Lambda
}}\sum
_{\theta
}2^{s}%
\mathbb{E}_{\mathbf{X}|\theta}d^{s}
\bigl( T,\psi( \theta) \bigr)
\nonumber
\\[-9pt]
\\[-9pt]
\nonumber
& =&\frac{1}{2^{r}D_{\Lambda}}\sum
_{\theta}\mathbb{E}_{\mathbf{X}%
|\theta} \bigl[ 2d \bigl( T,\psi(
\theta) \bigr) \bigr]^{s}.
\end{eqnarray}
Set $\hat{\theta}=\arg\min_{\theta\in\Theta}d^{s} ( T,\psi(\theta))$.
Note that the minimum is not necessarily unique.
When it is not unique, pick $\hat{\theta}$ to be any point in the minimum
set. Then the triangle inequality for the metric $d$ gives%
%
%
\begin{eqnarray}\label{Tri}
\mathbb{E}_{\mathbf{X}|\theta}d^{s} \bigl( \psi( \hat{\theta} )
,\psi(
\theta) \bigr) &\leq&\mathbb{E}_{\mathbf
{X}|\theta}%
\bigl[ d \bigl( \psi(
\hat{\theta} ) ,T \bigr) +d \bigl( T,\psi( \theta) \bigr) \bigr]^{s}
\nonumber
\\[-8pt]
\\[-8pt]
\nonumber
&\leq&\mathbb{E}_{\mathbf
{X}|\theta
} \bigl[ 2d \bigl( T,\psi( \theta) \bigr)
\bigr]^{s},
\end{eqnarray}
where the\vspace*{1pt} last inequality is due to the fact $d ( \psi(
\hat{\theta%
} ) ,T ) =d ( T,\psi( \hat{\theta} )
) \leq
d ( T,\psi( \theta) ) $ from the definition
of $\hat{%
\theta}$. Equations (\ref{MaxAve}) and (\ref{Tri}) together yield%
%
%
\begin{eqnarray}\label{part1}
\max_{\Theta}2^{s}\mathbb{E}_{\mathbf{X}|\theta}d^{s}
\bigl( T,\psi( \theta) \bigr)
&\geq&\frac{1}{2^{r}D_{\Lambda
}}\sum
_{\theta}%
\mathbb{E}_{\mathbf{X}|\theta}d^{s}
\bigl( \psi( \hat{\theta} ) ,\psi( \theta) \bigr)\nonumber\\
&\geq&\frac{1}{2^{r}D_{\Lambda}}\sum_{\theta}
\mathbb{E}_{\mathbf
{X}%
|\theta}\frac{d^{s} ( \psi( \hat{\theta} ) ,\psi
(\theta) ) }{H ( \gamma( \hat{\theta}) ,\gamma
( \theta) ) \vee1}\cdot H \bigl( \gamma( \hat{
\theta}) ,\gamma( \theta) \bigr)\\
&\geq&\alpha\cdot\frac{1}{2^{r}D_{\Lambda}}\sum_{\theta}\mathbb{E}_%
{\mathbf{X}|\theta}H \bigl( \gamma( \hat{\theta} ) ,\gamma(\theta) \bigr),\nonumber
\end{eqnarray}
where the last step follows from the definition of $\alpha$ in
equation (%
\ref{alpha}).

We now show%
%
%
\begin{equation}
\frac{1}{2^{r}D_{\Lambda}}\sum_{\theta}\mathbb{E}_{\mathbf
{X}|\theta
}H
\bigl( \gamma( \hat{\theta} ) ,\gamma( \theta) \bigr) \geq\frac{r}{2}
\min_{i}\llVert\bar\mathbb{P}_{i,0}\wedge\bar\mathbb{P}_{i,1}\rrVert,\label{part2}
\end{equation}
which immediately implies $\max_{\Theta}2^{s}\mathbb{E}_{\mathbf
{X}|\theta
}d^{s} ( T,\psi( \theta) ) \geq\alpha\frac
{r}{2}%
\min_{i}\llVert\bar\mathbb{P}_{i,0}\wedge\bar\mathbb{P}%
_{i,1}\rrVert,$ and Lemma~\ref{AL} follows. From the definition
of $H$
in equation (\ref{H}) we write
\[
\frac{1}{2^{r}D_{\Lambda}}\sum_{\theta}\mathbb{E}_{\mathbf
{X}|\theta
}H
\bigl( \gamma( \hat{\theta} ) ,\gamma( \theta) \bigr) =\frac
{1}{2^{r}D_{\Lambda}}\sum
_{\theta}\sum_{i=1}^{r}
\mathbb{E}_{%
\mathbf{X}|\theta}\bigl\vert\gamma_{i} ( \hat{\theta} ) -
\gamma_{i} ( \theta) \bigr\vert.
\]
The right-hand side can be further written as
\begin{eqnarray*}
&&\sum_{i=1}^{r}\frac{1}{2^{r}D_{\Lambda}}\sum
_{\rho\in\Gamma
} \biggl[ \sum_{ \{ \theta:\gamma(\theta)=\rho\} }
\mathbb{E}_{\mathbf{X}%
|\theta}\bigl\vert\gamma_{i}(\hat{\theta})-
\gamma_{i}(\theta)\bigr\vert\biggr]
\\
&&\qquad=\frac{1}{2}\sum_{i=1}^{r} \biggl[
\frac{1}{2^{r-1}D_{\Lambda}}%
\sum_{ \{ \rho:\rho_{i}=0 \} }\sum
_{ \{ \theta
:\gamma
(\theta)=\rho\} }\int\gamma_{i}(\hat{\theta})\,d\mathbb
{P}_{\theta}\\
&&\hspace*{61pt}{}+%
\frac{1}{2^{r-1}D_{\Lambda}}\sum
_{ \{ \rho:\rho_{i}=1 \}
}\sum_{ \{ \theta:\gamma(\theta)=\rho\} }\int\bigl(1-
\gamma_{i}(%
\hat{\theta})\bigr)\,d\mathbb{P}_{\theta^{\prime}}
\biggr]
\\
&&\qquad=\frac{1}{2}\sum_{i=1}^{r} \biggl[
\int\gamma_{i}(\hat{\theta}) \biggl(\frac{1}{%
2^{r-1}D_{\Lambda}}\sum
_{ \{ \rho:\rho_{i}=0 \} }\sum_{ \{
\theta:\gamma(\theta)=\rho\} }\,d
\mathbb{P}_{\theta}\biggr)\\
&&\hspace*{61pt}{}+\int\bigl(1-\gamma_{i}(\hat{\theta})
\bigr) \biggl(\frac{1}{2^{r-1}D_{\Lambda}}\sum_{ \{
\rho:\rho_{i}=1 \} }\sum
_{ \{ \theta:\gamma(\theta
)=\rho
\} }\,d\mathbb{P}_{\theta}\biggr) \biggr]
\\
&&\qquad=\frac{1}{2}\sum_{i=1}^{r} \biggl[
\int\gamma_{i}(\hat{\theta})\,d\bar\mathbb{P}_{i,0}+
\int\bigl(1-\gamma_{i}(\hat{\theta})\bigr)\,d
\bar\mathbb{P}_{i,1}%
\biggr] .
\end{eqnarray*}

The following elementary result is useful to establish the lower bound for
the minimax risk. See, for example, page 40 of \citet{LeCam73}.

%
\begin{lemma}
\label{aff} The total variation affinity satisfies
\[
\llVert\mathbb{P}\wedge\mathbb{Q}\rrVert=\inf_{0\leq
f\leq
1} \biggl\{ \int fd
\mathbb{P}+\int(1-f)\,d\mathbb{Q} \biggr\} .
\]
\end{lemma}

It follows immediately from Lemma~\ref{aff} that
\begin{eqnarray*}
\frac{1}{2}\sum_{i=1}^{r} \biggl[
\int\gamma_{i}(\hat{\theta})\,\bar\mathbb{P}_{i,0}+\int\bigl(1-\gamma_{i}(\hat{\theta})\bigr)\,d\bar\mathbb{P}_{i,1}
\biggr] &\geq&\frac{1}{2}\sum
_{i=1}^{r}\llVert\bar\mathbb{P}_{i,0}
\wedge\bar\mathbb{P}_{i,1}\rrVert\\
&\geq&\frac{r}{2}
\min_{i}\llVert\bar\mathbb{P}_{i,0}\wedge
\bar\mathbb{P}_{i,1}\rrVert,
\end{eqnarray*}
and so equation (\ref{part2}) is established.

\subsection{\texorpdfstring{Proof of Lemma \protect\ref{dffbd}}{Proof of Lemma 5}}
Let $v= ( v_{i} ) $ be a column $p$-vector with $v_{i}=0$
for $1\le
i \le p-r$ and $v_{i}=1$ for $p-r+1\leq i\leq p$, that is, $v= (
1 \{
p-r+1\leq i\leq p \} )_{p\times1}$. Set $w= (
w_{i} ) =%
[ \Sigma( \theta) -\Sigma( \theta^{\prime
} ) %
] v$. Note that for each $i$, if $\vert\gamma_{i}(\theta
)-\gamma_{i}(\theta^{\prime})\vert=1$, we have $\vert
w_{i}\vert=k\epsilon_{n,p}$. Then there are at least $H(\gamma
(\theta),\gamma(\theta^{\prime}))$ number of elements $w_{i}$ with $
\vert w_{i}\vert=k\epsilon_{n,p}$, which implies%
\[
\bigl\llVert\bigl[ \Sigma(\theta)-\Sigma\bigl(\theta^{\prime
}\bigr)
\bigr] v\bigr\rrVert_{2}^{2}\geq H\bigl(\gamma(\theta),
\gamma\bigl(\theta^{\prime
}\bigr)\bigr)\cdot( k\epsilon_{n,p}
)^{2}.
\]
Since $\llVert v\rrVert^{2}=r\leq p$, the equation above yields
\[
\bigl|\!\bigl|\!\bigl|\Sigma(\theta)-\Sigma\bigl(\theta^{\prime
}\bigr)
\bigr|\!\bigr|\!\bigr|^{2}\geq\frac{\llVert[ \Sigma(\theta
)-\Sigma(\theta
^{\prime
}) ]
v\rrVert_{2}^{2}}{\llVert v\rrVert^{2}}\geq\frac
{H(\gamma
(\theta),\gamma(\theta^{\prime}))\cdot(k\epsilon_{n,p})^{2}}{p},
\]
that is,%
\[
\frac{ |\!|\!|\Sigma(\theta)-\Sigma(\theta^{\prime})
|\!|\!|^{2}%
}{H(\gamma(\theta),\gamma(\theta^{\prime}))}\geq\frac{(k\epsilon
_{n,p})^{2}}{p}
\]
when $H(\gamma(\theta),\gamma(\theta^{\prime}))\geq1$. %

\subsection{\texorpdfstring{Proof of Lemma \protect\ref{affbd}}{Proof of Lemma 6}}

The proof of the bound for the affinity given in Lemma~\ref{affbd} is
involved. We break the proof into a few major technical lemmas which are
proved in Section~\ref{chisquaresec} and the supplementary material. Without loss of
generality we consider only the case $i=1$ and prove that there exists a
constant $c_{1}>0$ such that $\llVert\bar\mathbb{P}_{1,0}\wedge
\bar\mathbb{P}_{1,1}\rrVert\geq c_{1}$. The following lemma
is the
key step which turns the problem of bounding the total variation affinity
into a chi-squared distance calculation on Gaussian mixtures.%
%

%
\begin{lemma}
\label{chisquareslem} \textup{(i)} There exists a constant
$c_{2}<1$ such
that
%
%
\begin{equation}
\tilde\mathbb{E}_{ ( \gamma_{-1},\lambda_{-1} ) } \biggl\{ \int
\biggl( \frac{d\bar\mathbb{P}_{ ( 1,1,\gamma_{-1},\lambda
_{-1} )
}}{d\bar\mathbb{P}_{ ( 1,0,\gamma_{-1},\lambda_{-1} )
}}
\biggr)^{2}\,d\bar\mathbb{P}_{ ( 1,0,\gamma_{-1},\lambda_{-1} )
}-1 \biggr\} \leq
c_{2}^{2}. \label{chi-squarebd}
\end{equation}

\textup{(ii)} Moreover, equation \eqref{chi-squarebd} implies that $%
\llVert\bar\mathbb{P}_{1,0}\wedge\bar\mathbb{P}_{1,1}\rrVert
\geq1-c_{2}>0$.
\end{lemma}

The proof of Lemma~\ref{chisquareslem}(ii) is relatively easy and is given
in the supplementary material. Our goal in the remainder of this proof is to
establish %
\eqref{chi-squarebd}, which requires detailed understanding of
$\bar\mathbb{P}_{ ( 1,0,\gamma_{-1},\lambda_{-1} ) }$ and the mixture
distribution $\bar\mathbb{P}_{(1,1,\gamma_{-1},\lambda_{-1}) }$
as well
as a careful analysis of the cross-product terms in the chi-squared
distances on the left-hand side of \eqref{chi-squarebd}.\vadjust{\goodbreak}

From the definition of $\theta$ in equation (\ref{theta}) and
$\bar\mathbb{P}_{ ( 1,0,\gamma_{-1},\lambda_{-1} ) }$ in equation
(\ref%
{avepibd}), $\gamma_{1}=0$ implies $\bar\mathbb{P}_{ (
1,0,\gamma
_{-1},\lambda_{-1} ) }$ is a single multivariate normal distribution
with a covariance matrix,
%
%
\begin{equation}
\Sigma_{0}=\pmatrix{ 1 &
\mathbf{0}_{1\times( p-1 ) }
\vspace*{2pt}\cr
\mathbf{0}_{ ( p-1 ) \times1} & \mathbf{S}_{ (
p-1 )
\times( p-1 ) }}. \label{sigma0}
\end{equation}
Here $\mathbf{S}_{ ( p-1 ) \times( p-1 )
}= (
s_{ij} )_{2\leq i,j\leq p}$ is a symmetric matrix uniquely
determined by $ ( \gamma_{-1},\lambda_{-1} ) = ( (\gamma_{2},\ldots
,\gamma_{r}),(\lambda_{2},\ldots,\lambda_{r}) ) $ where for $i\le j$,
\[
s_{ij}=\cases{ %
1, &\quad$ i=j$,
\vspace*{2pt}\cr
\epsilon_{n,p}, & \quad$\gamma_{i}= \lambda_{i}
( j ) =1$,
\vspace*{2pt}\cr
0, & \quad $\mbox{otherwise}.$}
\]
Let %
\[
\Lambda_{1} ( c ) = \bigl\{ a
\in B\dvtx\exists\theta\in\Theta\mbox{ such that }\lambda_{1}(
\theta)=a\mbox{ and }\lambda_{-1}(\theta)=c \bigr\} ,
\]
which gives the set of all possible values of the first row with the
rest of
the rows fixed, that is, $\lambda_{-1}(\theta)=c$. Let $n_{\lambda
_{-1}}$ be
the number of columns of $\lambda_{-1}$ with the column sum equal to $2k$
for which the first row has no choice but to take value $0$ in this column.
Set $p_{\lambda_{-1}}=r-n_{\lambda_{-1}}$. It is helpful to observe
that $%
p_{\lambda_{-1}}\geq p/4-1$. Since $n_{\lambda_{-1}}\cdot2k\leq
r\cdot k$%
, the total number of $1$s in the upper triangular matrix by the
construction of the parameter set, we thus have $n_{\lambda_{-1}}\leq r/2$,
which immediately implies $p_{\lambda_{-1}}=r-n_{\lambda_{-1}}\geq
r/2\geq
p/4-1$. It follows $\operatorname{Card}( \Lambda_{1} ( \lambda_{-1}
) ) =({{p_{\lambda_{-1}}\atop k}})$. Then, from the
definitions in equations (\ref{theta}) and (\ref{avepibd}), $\bar\mathbb{P}
_{ ( 1,1,\gamma_{-1},\lambda_{-1} ) }$ is an average of $({%
p_{\lambda_{-1}}\atop k})$ multivariate normal distributions with covariance
matrices of the following form:%
%
%
\begin{equation}
\pmatrix{ 1 & \mathbf{r}_{1\times( p-1 ) }
\vspace*{2pt}\cr
( \mathbf{r}_{1\times( p-1 ) } )^{T} & \mathbf{S}%
_{ ( p-1 ) \times( p-1 ) }} ,\label{matrixform}
\end{equation}
where $\llVert\mathbf{r}\rrVert_{0}=k$ with nonzero
elements of $r$
equal $\epsilon_{n,p}$ and the submatrix $\mathbf{S}_{(p-1)\times(p-1)}$
is the same as the one for $\Sigma_{0}$ given in (\ref{sigma0}).

Recall that for each $\theta\in\Theta$, $\mathbb{P}_\theta$ is the joint
distribution of the $n$ i.i.d. multivariate normal variables $\mathbf
{X}_{1},\ldots, \mathbf{X}_{n}$. So each term in the chi-squared
distance on
the left-hand side of \eqref{chi-squarebd} is of the form $ (\int
\frac{%
g_{1}g_{2}}{g_{0}} )^n$ where $g_{i}$ are the density function of
$%
N (0,\Sigma_{i} ) $ for $i=0,1$ and $2$, with $\Sigma_0$ defined
in \eqref{sigma0} and $\Sigma_{1}$ and $\Sigma_{2}$ of the form~(\ref{matrixform}).

The following lemma is useful for calculating the cross product terms
in the
chi-squared distance between Gaussian mixtures. The proof of the lemma is
straightforward and is thus omitted.

%
\begin{lemma}
\label{crossproduct} Let $g_{i}$ be the density function of $N (
0,\Sigma_{i} ) $ for $i=0,1$ and $2$, respectively. Then
\[
\int\frac{g_{1}g_{2}}{g_{0}}= \bigl[ \det\bigl( I-\Sigma_{0}^{-2}
( \Sigma_{1}-\Sigma_{0} ) ( \Sigma_{2}-
\Sigma_{0} ) \bigr) \bigr]^{-{{1}/{2}}}.
\]
\end{lemma}

Let $\Sigma_{0}$ be defined in \eqref{sigma0} and determined by
$ (
\gamma_{-1},\lambda_{-1} ) $. Let $\Sigma_{1}$ and $\Sigma_{2}$ be
of the form (\ref{matrixform}) with the first row $\lambda_{1}$ and $
\lambda_{1}^{\prime}$, respectively. Set
%
%
\begin{equation}
R_{\lambda_{1},\lambda_{1}^{\prime}}^{\gamma_{-1},\lambda
_{-1}}=-\log\det\bigl( I-\Sigma_{0}^{-2}
( \Sigma_{0}-\Sigma_{1} ) ( \Sigma_{0}-
\Sigma_{2} ) \bigr) . \label{R}
\end{equation}
We sometimes drop the indices $(\lambda_{1}$, $\lambda_{1}^{\prime
})$ and
$ ( \gamma_{-1},\lambda_{-1} ) $ from $\Sigma_{i}$ to simplify
the notation whenever there is no ambiguity. Then each term in the
chi-squared distance on the left-hand side of \eqref{chi-squarebd} can be
expressed as in the form of
\[
\exp\biggl( \frac{n}{2}\cdot R_{\lambda_{1},\lambda_{1}^{\prime
}}^{\gamma
_{-1},\lambda_{-1}} \biggr)
-1.
\]

Define
\begin{eqnarray*}
&&\Theta_{-1} ( a_{1},a_{2} ) = \{ 0,1
\}^{r-1}\otimes\bigl\{ c\in\Lambda_{-1}\dvtx\exists
\theta_{i}\in\Theta,
i=1,2,
\\
&&\hspace*{125pt}\mbox{such that }\lambda_{1}(\theta_{i})=a_{i},
\lambda_{-1}(\theta_{i})=c \bigr\} .
\end{eqnarray*}
It is a subset of $\Theta_{-1}$ in which the element can pick both $a_{1}$
and $a_{2}$ as the first row to form parameters in $\Theta$. From
Lemma \ref%
{crossproduct} the average of the chi-squared distance on the left-hand side
of equation (\ref{chi-squarebd}) can now be written as
%
%
\begin{eqnarray}\label{chisquare}
&&\tilde\mathbb{E}_{ ( \gamma_{-1},\lambda_{-1} )
} \biggl\{ \tilde\mathbb{E}_{ ( \lambda_{1},\lambda_{1}^{\prime
} )
|\lambda_{-1}}
\biggl[ \exp\biggl(\frac{n}{2}\cdot R_{\lambda_{1},\lambda
_{1}^{\prime}}^{\gamma_{-1},\lambda_{-1}}
\biggr)-1 \biggr] \biggr\}
\nonumber
\\[-8pt]
\\[-8pt]
\nonumber
&&\qquad=\tilde\mathbb{E}_{ ( \lambda_{1},\lambda_{1}^{\prime
} )
} \biggl\{ \tilde\mathbb{E}_{ ( \gamma_{-1},\lambda
_{-1} )
| ( \lambda_{1},\lambda_{1}^{\prime} ) }
\biggl[ \exp\biggl(\frac{n}{2}%
\cdot R_{\lambda_{1},\lambda_{1}^{\prime}}^{\gamma_{-1},\lambda
_{-1}}
\biggr)-1%
\biggr] \biggr\},
\end{eqnarray}
where $\lambda_{1}$ and $\lambda_{1}^{\prime}$ are independent and
uniformly distributed over $\Lambda_{1} ( \lambda_{-1} ) $ (not
over $B$) for given $\lambda_{-1}$, and the distribution of $ (
\gamma_{-1},\lambda_{-1} ) $ given $ ( \lambda_{1},\lambda
_{1}^{\prime
} ) $ is uniform over $\Theta_{-1}$ $ ( \lambda_{1},\lambda
_{1}^{\prime} ) $, but the marginal distribution of $\lambda_{1}$ and
$\lambda_{1}^{\prime}$ are not independent and uniformly distributed
over $%
B$.

Let $\Sigma_{1}$ and $\Sigma_{2}$ be two covariance matrices of the
form (\ref{matrixform}). Note that $\Sigma_{1}$ and $\Sigma_{2}$ differ from
each other only in the first row/column. Then $\Sigma_{i}-\Sigma
_{0}$, $%
i=1 $ or $2$, has a very simple structure. The nonzero elements only appear
in the first row/column, and in total there are at most $2k$ nonzero
elements. This property immediately implies the following lemma which makes
the problem of studying the determinant in Lemma~\ref{crossproduct}
relatively easy. The proof of Lemma~\ref{sigmai} below is given in the
supplementary material.

%
\begin{lemma}
\label{sigmai} Let $\Sigma_{0}$ be defined in \eqref{sigma0} and let
$%
\Sigma_{1}$ and $\Sigma_{2}$ be two covariance matrices of the form
(\ref{matrixform}). Define $J$ to be the number of overlapping $\epsilon
_{n,p}$%
's between $\Sigma_{1}$ and $\Sigma_{2}$ on the first row, and
\[
Q\stackrel{\bigtriangleup} {=} ( q_{ij} )_{1\leq i,j\leq
p}= (
\Sigma_{1}-\Sigma_{0} ) ( \Sigma_{2}-
\Sigma_{0} ) .
\]
There are index subsets $I_{r}$ and $I_{c}$ in $ \{ 2,\ldots
,p \} $
with $\operatorname{Card} ( I_{r} ) =\operatorname{Card} (
I_{c} ) =k$
and $\operatorname{Card} ( I_{r}\cap I_{c} ) =J$ such that
\[
q_{ij}=\cases{ %
J\epsilon_{n,p}^{2},
&\quad $i=j=1$,
\vspace*{2pt}\cr
\epsilon_{n,p}^{2}, &\quad $i\in I_{r}\mbox{ and }j\in
I_{c}$,
\vspace*{2pt}\cr
0, & \quad$\mbox{otherwise},$}
\]
and the matrix $ ( \Sigma_{0}-
\Sigma_{1} ) ( \Sigma_{0}-\Sigma_{2} ) $ has rank
$2$ with two identical nonzero eigenvalues $J\epsilon_{n,p}^{2}$
when $J>0$.
\end{lemma}
The matrix $Q$ is determined by two interesting parts,
the first element $%
q_{11}=J\epsilon_{n,p}^{2}$
and a very special $k\times k$ square matrix $%
( q_{ij}\dvtx
i\in I_{r}\mbox{ and }j\in I_{c} ) $ with all elements
equal to $\epsilon_{n,p}^{2}$. The following result, which is
proved in the supplementary material, shows that $R_{\lambda_{1},\lambda
_{1}^{\prime}}^{\gamma
_{-1},\lambda_{-1}}$ is
approximately equal to %
\[
-\log\det\bigl( I- ( \Sigma_{0}-
\Sigma_{1} ) ( \Sigma_{0}-\Sigma_{2} ) \bigr) =-2
\log\bigl( 1-J\epsilon_{n,p}^{2} \bigr),
\]
where $J$ is defined in Lemma~\ref{sigmai}. Define%
\begin{eqnarray*}
\Lambda_{1,J}&= &\bigl\{ \bigl( \lambda_{1},
\lambda_{1}' \bigr) \in B\otimes B\dvtx\mbox{the number
of overlapping }\epsilon_{n,p}\mbox{'s between}\\
&&\hspace*{208pt}\lambda_{1}\mbox{ and }\lambda_{1}^{\prime}\mbox{ is }J
\bigr\} .
\end{eqnarray*}

%
\begin{lemma}
\label{Rlem}
Let $R_{\lambda_{1},\lambda_{1}^{\prime}}^{\gamma_{-1},\lambda_{-1}}$
be defined in equation (\ref{R}). Then%
%
%
\begin{equation}
R_{\lambda_{1},\lambda_{1}^{\prime}}^{\gamma_{-1},\lambda
_{-1}}=-2\log\bigl( 1-J\epsilon_{n,p}^{2}
\bigr) +R_{1,\lambda_{1},\lambda
_{1}^{\prime
}}^{\gamma_{-1},\lambda_{-1}}, \label{Rdecomp}
\end{equation}
where $R_{1,\lambda_{1},\lambda_{1}^{\prime}}^{\gamma_{-1},\lambda
_{-1}} $ satisfies, uniformly over all $J$,
%
%
\begin{equation}
\tilde\mathbb{E}_{ ( \lambda_{1},\lambda_{1}^{\prime}
) |J}%
\biggl[ \tilde\mathbb{E}_{ ( \gamma_{-1},\lambda_{-1}
) | (
\lambda_{1},\lambda_{1}^{\prime} ) }\exp\biggl( \frac{n}{2}%
R_{1,\lambda_{1},\lambda_{1}^{\prime}}^{\gamma_{-1},\lambda
_{-1}} \biggr) \biggr] \leq{\frac{3}{2}}.
\label{R1eq}
\end{equation}
\end{lemma}

With the preparations given above, we are now ready to establish
equation (\ref{chi-squarebd}) and thus complete the proof of Lemma
\ref{affbd}.

\subsubsection*{Proof of equation (\protect\ref{chi-squarebd})}
\label{chisquaresec}

Equation (\ref{Rdecomp}) in Lemma~\ref{Rlem} yields that
\begin{eqnarray*}
&&\tilde\mathbb{E}_{ ( \lambda_{1},\lambda_{1}^{\prime
} )
} \biggl\{ \tilde\mathbb{E}_{ ( \gamma_{-1},\lambda
_{-1} )
| ( \lambda_{1},\lambda_{1}^{\prime} ) }
\biggl[ \exp\biggl(\frac{n}{2}%
R_{\lambda_{1},\lambda_{1}^{\prime}}^{\gamma_{-1},\lambda_{-1}}
\biggr)-1%
\biggr] \biggr\}
\\
&&\qquad=\tilde\mathbb{E}_{J} \biggl\{ \exp\bigl[ -n\log\bigl
( 1-J
\epsilon_{n,p}^{2} \bigr) \bigr] \\
&&\hspace*{50pt}{}\times\tilde\mathbb{E}_{ ( \lambda
_{1},\lambda_{1}^{\prime} ) |J} \biggl[ \tilde\mathbb{E}_{ (
\gamma_{-1},\lambda_{-1} ) | ( \lambda_{1},\lambda
_{1}^{\prime
} ) }\exp\biggl(
\frac{n}{2}R_{1,\lambda_{1},\lambda
_{1}^{\prime
}}^{\gamma_{-1},\lambda_{-1}} \biggr) \biggr] -1 \biggr\} .
\end{eqnarray*}

Recall that $J$ is the number of overlapping $\epsilon_{n,p}$'s
between $%
\Sigma_{1}$ and $\Sigma_{2}$ on the first row. It is easy to see that $J$
has the hypergeometric distribution as $\lambda_{1}$ and $\lambda
_{1}^{\prime}$ vary in $B$ for each given $\lambda_{-1}$. For $0\leq
j\leq
k$,
%
%
\begin{eqnarray}\label{plambda}
\tilde\mathbb{E}_{J}\bigl(\mathbf{1} \{ J=j \} |
\lambda_{-1}\bigr) &=&%
\pmatrix{k
\cr
j}\pmatrix{p_{\lambda_{-1}}-k
\cr
k-j}\Big/\pmatrix{p_{\lambda_{-1}}
\cr
k}
\nonumber
\\
&=&\frac{ ( {k!}/{ ( k-j ) !} )^{2}}{{(p_{\lambda
_{-1}}! ( p_{\lambda_{-1}}-2k+j ) !)}/{ [ (
p_{\lambda
_{-1}}-k ) ! ]^{2}}}\cdot\frac{1}{j!}\\
&\leq&\biggl( \frac
{k^{2}}{%
p_{\lambda_{-1}}-k}
\biggr)^{j},
\nonumber
\end{eqnarray}
where $\frac{k!}{ ( k-j ) !}$ is a product of $j$ term with each
term $\leq k$ and for $\frac{p_{\lambda_{-1}}! ( p_{\lambda
_{-1}}-2k+j ) !}{ [ ( p_{\lambda_{-1}}-k )
! ]^{2}}
$ it is bounded below by a product of $j$ term with each term $\geq
p_{\lambda_{-1}}-j$. Since $p_{\lambda_{-1}}\geq p/4-1$ for all
$\lambda_{-1}$, we have
\[
\tilde\mathbb{E}\bigl(\mathbf{1} \{ J=j \} \bigr)=\tilde\mathbb{E}%
_{\lambda_{-1}} \bigl[ \tilde\mathbb{E}_{J} \bigl(
\mathbf{1} \{ J=j \} |\lambda_{-1} \bigr) \bigr] \leq\biggl(
\frac
{k^{2}}{p/4-1-k}%
\biggr)^{j}.
\]
Thus%
%
%
\begin{eqnarray} \label{limitkp}
&&\tilde\mathbb{E}_{ ( \gamma_{-1},\lambda_{-1} )
} \biggl\{ \int\biggl( \frac{d\bar\mathbb{P}_{ ( 1,1,\gamma
_{-1},\lambda
_{-1} )
}}{d\bar\mathbb{P}_{ ( 1,0,\gamma_{-1},\lambda_{-1} )
}}
\biggr)^{2}\,d\bar\mathbb{P}_{ ( 1,0,\gamma_{-1},\lambda_{-1} )
}-1 \biggr\}
\nonumber
\\
&&\qquad\leq\sum_{j\geq0} \biggl( \frac{k^{2}}{p/4-1-k}
\biggr)^{j} \biggl\{ \exp%
\bigl[ -n\log\bigl( 1-j
\epsilon_{n,p}^{2} \bigr) \bigr] \cdot{\frac{3}{2}%
}-1 \biggr\}
\\
&&\qquad={\frac{3}{2}}\sum_{j\geq1} \biggl(
\frac{k^{2}}{p/4-1-k} \biggr)^{j}\exp%
\bigl[ 2j \bigl(
\upsilon^{2}\log p \bigr) \bigr]
\nonumber\\
&&\qquad\quad{}+ \biggl( \frac{k^{2}}{p/4-1-k} \biggr)^{0} \biggl
\{ \exp\bigl[ -n
\log\bigl( 1-0\cdot\epsilon_{n,p}^{2} \bigr) \bigr] \cdot{
\frac
{3}{2}}-1 \biggr\}
\nonumber
\\
&&\qquad\leq C\sum_{j\geq1} \bigl( p^{{(\beta-1)}/{\beta
}}\cdot
p^{-2\upsilon
^{2}} \bigr)^{-j}+{\frac{1}{2}}<C\sum
_{j\geq1} \bigl( p^{{(\beta-1)}/{(2\beta)}} \bigr)^{-j}+
\frac{1}{2}<c_{2}^{2}
\nonumber
\end{eqnarray}
by setting $c_{2}^{2}=3/4$, where the last step follows from $\upsilon
^{2}<%
\frac{\beta-1}{54\beta}$ and $k^{2}=O ( \frac{n}{\log p} )
=O ( \frac{p^{1/\beta}}{\log p} ) $ as defined in Section
\ref{lowbdsec}.

%
\begin{remark}
The condition $p\geq n^{\beta}$ for some $\beta>1$ is
assumed so
that%
\[
\frac{k^{2}}{p_{\lambda_{-1}}-k}\leq\frac{k^{2}}{p/4-k}=\frac
{O (
n/\log p ) }{p/4-k}=o \bigl(
p^{-\varepsilon} \bigr)
\]
for some $\varepsilon>0$ to make the term (\ref{limitkp}) to be
$o (
1 ) $.
\end{remark}

\subsection{\texorpdfstring{Proof of Theorem \protect\ref{MinimaxOpe1}}{Proof of Theorem 3}}
\label{uppersec}

The following lemma, which is proved in Cai and Zhou (\citeyear{CaiZho}), is now
useful to
prove Theorem~\ref{MinimaxOpe1}.

%
\begin{lemma}
\label{A0bound} Define the event $A_{ij}$ by
%
%
\begin{equation}
A_{ij}= \biggl\{ \vert\hat{\sigma}_{ij}-
\sigma_{ij}\vert\leq4\min\biggl\{ \vert\sigma_{ij}
\vert,\gamma\sqrt{\frac{\log p}{%
n}} \biggr\} \biggr\} . \label{Aij}
\end{equation}
Then $\mathbb{P} ( A_{ij} ) \geq1-2C_{1}p^{-9/2}.$
\end{lemma}

Let $D= ( d_{ij} )_{1\leq i,j\leq p}$ with $d_{ij}= (
\hat{%
\sigma}_{ij}-\sigma_{ij} ) I(A_{ij}^{c})$. Then
%
%
\begin{eqnarray}\label{thm31bnd}
&&\mathbb{E}_{\mathbf{X}|\theta} |\!|\!|\hat{\Sigma}-\Sigma|\!|\!
|^{2} \nonumber\\
&&\qquad
\leq2\mathbb{E}_{\mathbf{X}|\theta} \biggl[ \sup_{j}\sum
_{i\neq
j}\vert\hat{\sigma}_{ij}-
\sigma_{ij}\vert I(A_{ij}) \biggr]^{2}+2
\mathbb{E}_{\mathbf{X}|\theta} |\!|\!|D |\!|\!|_{1}^{2}+C
\frac{%
\log p}{n}
\\
&&\qquad\leq 32 \biggl[ \sup_{j}\sum_{i\neq j}\min
\biggl\{ \vert\sigma_{ij}\vert,\gamma\sqrt{\frac{\log p}{n}}
\biggr\} \biggr]^{2}+2%
\mathbb{E}_{\mathbf{X}|\theta} |\!|\!|D
|\!|\!|_{1}^{2}+C\frac{\log p%
}{n}.\nonumber
\end{eqnarray}

We will see that the first term in equation (\ref{thm31bnd}) is dominating
and is bounded by $Cc_{n,p}^{2} ( \frac{\log p}{n} )^{1-q}$, while
the second term $\mathbb{E}_{\mathbf{X}|\theta} |\!|\!|D
|\!|\!|_{1}^{2}$ is negligible.

Set $k^{\ast}=\lfloor c_{n,p} ( {\frac{n}{\log p}} )^{q/2}\rfloor
$. Then we have
%
%
\begin{eqnarray}\label{boundDominate}
\sum_{i\neq j}\min\biggl\{ \vert
\sigma_{ij}\vert,\gamma\sqrt{%
\frac{\log p}{n}} \biggr\}
&\leq&\gamma\biggl( \sum_{i\leq k^{\ast
}}+\sum
_{i>k^{\ast}} \biggr) \min\biggl\{ \vert\sigma_{
[ i ]
j}
\vert,\sqrt{\frac{\log p}{n}} \biggr\}
\nonumber
\\
&\leq&C_{5}k^{\ast}\sqrt{\frac{\log p}{n}}+C_{5}
\sum_{i>k^{\ast
}} \biggl( \frac{c_{n,p}}{i}
\biggr)^{{1}/{q}}
\nonumber
\\[-8pt]
\\[-8pt]
\nonumber
&\leq&C_{6} \biggl[ k^{\ast}\sqrt{\frac{\log p}{n}}+c_{n,p}^{{1}/{q}%
}
\cdot\bigl( k^{\ast} \bigr)^{1-{{1}/{q}}} \biggr]\\
& \leq&
C_{7}c_{n,p} \biggl( \frac{\log p}{n} \biggr)^{{(1-q)}/{2}},
\nonumber
\end{eqnarray}
which immediately implies equation (\ref{rateOper2}) if $\mathbb
{E}_{\mathbf{X}|\theta} |\!|\!|D |\!|\!|_{1}^{2}=O
( \frac{1}{n} ) $. We
shall now show that $\mathbb{E}_{\mathbf{X}|\theta} |\!|\!|
D |\!|\!|_{1}^{2}=O ( \frac{1}{n} ) $. Note that
\begin{eqnarray*}
\mathbb{E}_{\mathbf{X}|\theta} |\!|\!|D |\!|\!|_{1}^{2} &\leq
&p\sum_{ij}\mathbb{E}_{\mathbf{X}|\theta}d_{ij}^{2}\\
&=& p
\sum_{ij}\mathbb{E}_{%
\mathbf{X}|\theta} \bigl\{ \bigl[
d_{ij}^{2}I\bigl(A_{ij}^{c}\cap\bigl\{
\hat{%
\sigma}_{ij}=\sigma_{ij}^{\ast}
\bigr\} \bigr)+d_{ij}^{2}I\bigl(A_{ij}^{c}
\cap\{ \hat{\sigma}_{ij}=0 \}\bigr) \bigr] \bigr\}
\\
&=&p\sum_{ij}\mathbb{E}_{\mathbf{X}|\theta} \bigl\{
\bigl( \sigma_{ij}^{\ast}-\sigma_{ij}
\bigr)^{2}I\bigl(A_{ij}^{c}\bigr) \bigr\} +p\sum
_{ij}%
\mathbb{E}_{\mathbf{X}|\theta}
\sigma_{ij}^{2}I\bigl(A_{ij}^{c}\cap\{
\hat{%
\sigma}_{ij}=0 \}\bigr)
\\
&\equiv& R_{1}+R_{2}.
\end{eqnarray*}

Lemma~\ref{A0bound} yields
that $\mathbb{P} ( A_{ij}^{c} ) \leq2C_{1}p^{-9/2}$, and the
Whittle inequality implies $\sigma_{ij}^{\ast}-\sigma_{ij}$ has all
finite moments [cf. \citet{Whi60}] under the subgaussianity condition
(\ref{subGau}).
Hence
\begin{eqnarray*}
R_{1} &=&p\sum_{ij}\mathbb{E}_{\mathbf{X}|\theta}
\bigl\{ \bigl( \sigma_{ij}^{\ast}-\sigma_{ij}
\bigr)^{2}I\bigl(A_{ij}^{c}\bigr) \bigr\} \leq p
\sum_{ij}%
\bigl[ \mathbb{E}_{\mathbf{X}|\theta}
\bigl( \sigma_{ij}^{\ast
}-\sigma_{ij}
\bigr)^{6} \bigr]^{1/3}\mathbb{P}^{2/3} \bigl(
A_{ij}^{c} \bigr)
\\
&\leq&C_{8}p\cdot p^{2}\cdot\frac{1}{n}\cdot
p^{-3}=C_{8}/n.
\end{eqnarray*}
On the other hand,
\begin{eqnarray*}
R_{2} &=&p\sum_{ij}\mathbb{E}_{\mathbf{X}|\theta}
\sigma_{ij}^{2}I\bigl(A_{ij}^{c}\cap\{ \hat{
\sigma}_{ij}=0 \} \bigr)\\
&=& p\sum_{ij}%
\mathbb{E}_{\mathbf{X}|\theta}\sigma_{ij}^{2}I\biggl(\vert
\sigma_{ij}\vert\geq4\gamma\sqrt{\frac{\log p}{n}}\biggr)I\biggl(|
\sigma_{ij}^{\ast
}|\leq\gamma\sqrt{\frac{\log p}{n}}\biggr)
\\
&\leq&p\sum_{ij}\sigma_{ij}^{2}
\mathbb{E}_{\mathbf{X}|\theta
}I\biggl(\vert\sigma_{ij}\vert\geq4\gamma
\sqrt{\frac{\log p}{n}}\biggr)I\biggl(\vert\sigma_{ij}\vert
-\bigl
\vert\sigma_{ij}^{\ast}-\sigma_{ij}\bigr\vert
\leq\gamma\sqrt{\frac{\log p}{n}}\biggr)
\\
&\leq&p\sum_{ij}\sigma_{ij}^{2}
\mathbb{E}_{\mathbf{X}|\theta
}I\biggl(\bigl\vert\sigma_{ij}^{\ast}-
\sigma_{ij}\bigr\vert>\frac{3}{4}\vert\sigma_{ij}
\vert\biggr)I\biggl(\vert\sigma_{ij}\vert\geq4\gamma
\sqrt{%
\frac{\log p}{n}}\biggr)
\\
&\leq&\frac{p}{n}\sum_{ij}n
\sigma_{ij}^{2}C_{1}\exp\biggl( -{
\frac
{9}{%
2\gamma^{2}}}n\sigma_{ij}^{2} \biggr) I\biggl(\vert
\sigma_{ij}\vert\geq4\gamma\sqrt{\frac{\log p}{n}}\biggr)
\\
&=&\frac{p}{n}\sum_{ij} \biggl[ n
\sigma_{ij}^{2}\cdot C_{1}\exp\biggl(
-{%
\frac{1}{2\gamma^{2}}}n\sigma_{ij}^{2} \biggr)
\biggr] \cdot\exp\biggl( -{%
\frac{4}{\gamma^{2}}}n\sigma_{ij}^{2}
\biggr) I\biggl(\vert\sigma_{ij}\vert\geq4\gamma\sqrt{
\frac{\log p}{n}}\biggr)
\\
&\leq&C_{9}\frac{p}{n}\cdot p^{2}\cdot
p^{-16}\leq C_{9}/n.
\end{eqnarray*}
Putting $R_{1}$ and $R_{2}$ together yields that for some constant
$C>0$ ,
%
%
\begin{equation}
\mathbb{E}_{\mathbf{X}|\theta} |\!|\!|D |\!|\!|_{1}^{2}\leq{
\frac{C%
}{n}}. \label{D1}
\end{equation}
Theorem~\ref{MinimaxOpe1} is proved by combining equations (\ref
{thm31bnd}%
), (\ref{boundDominate}) and (\ref{D1}).

\subsection{\texorpdfstring{Proof of Theorem \protect\ref{MinimaxBD}}{Proof of Theorem 4}}
\label{Bregmanproofsec}

We establish separately the lower and upper bounds under the Bregman
divergence losses. The following lemma relates a general Bregman divergence
to the squared Frobenius norm.

%
\begin{lemma}
\label{FroBD}Assume that all eigenvalues of two symmetric matrices $X$
and $%
Y $ belong to $ [ \epsilon_{2},M_{2} ] $. Then there exist
constants $c_2>c_1>0$ depending only on $\epsilon_2$ and $M_2$ such
that for all $\phi\in
\Phi$ defined in (\ref{assumpphi}),
\[
c_1 |\!|\!|X-Y |\!|\!|_{F}^{2} \le
D_{\phi} ( X,Y ) \le c_2 |\!|\!|X-Y |\!|\!|_{F}^{2}.
\]
\end{lemma}

\begin{pf} Let the eigen
decompositions of $X$ and $Y$ be%
\[
X=\sum_{i=1}^{p}\lambda_{i}v_{i}^{T}v_{i}
\quad\mbox{and}\quad Y=\sum_{i=1}^{p}
\gamma_{i}u_{i}^{T}u_{i}.
\]
For every $\phi( X ) =\sum_{i=1}^{p}\varphi(
\lambda_{i} ) $ it is easy to see that
%
%
\begin{equation}
D_{\phi} ( X,Y ) =\sum_{i,j} \bigl(
v_{i}^{T}u_{i} \bigr)^{2} \bigl[
\varphi( \lambda_{i} ) -\varphi( \gamma_{j} ) -\varphi
{\acute{}}%
( \gamma_{j}
) \cdot( \lambda_{i}-\gamma_{j} ) %
\bigr]
.\label{TaylorBD}
\end{equation}
See \citet{KulSusDhi09}, Lemma 1. The Taylor expansion
gives%
\[
D_{\phi} ( X,Y ) =\sum_{i,j} \bigl(
v_{i}^{T}u_{i} \bigr)^{2}
\frac{%
1}{2}\varphi^{\prime\prime} ( \xi_{ij} ) (
\lambda_{i}-\gamma_{j} )^{2},
\]
where $\xi_{ij}$ is in between $\lambda_{i}$ and $\gamma_{j}$ and then
contained in $ [ \epsilon_{2},M_{2} ] $. From the
assumption in (%
\ref{assumpphi}), there are constants $c_{L}$ and $c_{u}$ such that $%
c_{L}\leq\varphi^{\prime\prime} ( \lambda) \leq c_{u}
$ for all
$\lambda$ in $ [ \epsilon_{2},M_{2} ] $, which immediately
implies
\begin{eqnarray*}
\frac{1}{2}c_L \sum_{i,j} \bigl(
v_{i}^{T}u_{i} \bigr)^{2} (
\lambda_{i}-\gamma_{j} )^{2}&\le& D_{\phi}
( X,Y ) \\
&\le&\frac{1}{2}c_u \sum_{i,j} \bigl(
v_{i}^{T}u_{i} \bigr)^{2} (
\lambda_{i}-\gamma_{j} )^{2}= |\!|\!|X-Y
|\!|\!|_{F}^{2}
\end{eqnarray*}
or equivalently
\[
\tfrac{1}{2}c_L |\!|\!|X-Y |\!|\!|_{F}^{2}\le
D_{\phi} ( X,Y ) \le\tfrac{1}{2}c_u |\!|\!|X-Y
|\!|\!|_{F}^{2}.
\]

\textit{Lower bound under Bregman matrix divergences}.
\label{secFrolowbd}
It is trivial to see that%
\[
\inf_{\hat{\Sigma}}\sup_{\theta\in\mathcal{P}_{q}(\tau
,c_{n,p})}\mathbb{E}_{\mathbf{X%
}|\theta}
\mathrm{L}_{\phi} \bigl( \hat{\Sigma},\Sigma( \theta) \bigr)
\geq c
\frac{1}{n}
\]
by constructing a parameter space with only diagonal matrices. It is then
enough to show that there exists some constant $c>0$ such that
\[
\inf_{\hat{\Sigma}}\max_{\mathcal{F}^{\ast}}\mathbb{E}_{\mathbf
{X}|\theta}%
\mathrm{L}_{\phi} \bigl( \hat{\Sigma},\Sigma( \theta) \bigr)
\geq
cc_{n,p} \biggl( \frac{\log p}{n} \biggr)^{1-q/2}
\]
for all $\phi\in\Phi$ defined in (\ref{assumpphi}). Equation (\ref%
{TaylorBD}) implies
%
%
\begin{equation}
\inf_{\hat{\Sigma}}\max_{\mathcal{F}^{\ast}}\mathbb{E}_{\mathbf
{X}|\theta}%
\mathrm{L}_{\phi} \bigl( \hat{\Sigma},\Sigma( \theta) \bigr) =
\inf_{\hat{\Sigma}\dvtx\epsilon_{1}I\prec\hat{\Sigma}\prec
2\tau
I}\max_{%
\mathcal{F}^{\ast}}\mathbb{E}_{\mathbf{X}|\theta}
\mathrm{L}_{\phi
} \bigl( \hat{\Sigma},\Sigma( \theta) \bigr) \label
{RestrictedMinimax}.
\end{equation}
Convexity of $\varphi$ implies $\varphi( \lambda_{i} )
-\varphi
( \gamma_{j} ) -\varphi^{\prime} ( \gamma_{j} ) \cdot
( \lambda_{i}-\gamma_{j} ) $ is nonnegative and
increasing when $%
\lambda_{i}$ moves away from the range $ [ \epsilon_{1},2\tau
] $
of those eigenvalues $\gamma_{j}$'s of $\Sigma( \theta) $.
From Lemma~\ref{FroBD} there is a universal constant $c_{L}$ such that%
\begin{eqnarray*}
&&\inf_{\hat{\Sigma}:\epsilon_{1}I\prec\hat{\Sigma}\prec2\tau
I}\max_{%
\mathcal{F}^{\ast}}\mathbb{E}_{\mathbf{X}|\theta}
\mathrm{L}_{\phi
} \bigl( \hat{\Sigma},\Sigma( \theta) \bigr)\\
&&\qquad
\geq c_{L}\inf_{\hat{%
\Sigma}:\epsilon_{1}I\prec\hat{\Sigma}\prec2\tau I}\max
_{\mathcal{F}%
^{\ast}}
\mathbb{E}_{\mathbf{X}|\theta}\frac{1}{p} \bigl|\!\bigl|\!\bigl|\hat{\Sigma}%
-
\Sigma( \theta) \bigr|\!\bigr|\!\bigr|_{F}^{2}
\\
&&\qquad=\frac{c_{L}}{p}\inf_{\hat{\Sigma}}\max_{\mathcal{F}^{\ast
}}
\mathbb{E}_{%
\mathbf{X}|\theta} \bigl|\!\bigl|\!\bigl|\hat{\Sigma}-\Sigma( \theta)
\bigr|\!\bigr|\!\bigr|_{F}^{2},
\end{eqnarray*}
where the last equality is from the same argument for equation (\ref%
{RestrictedMinimax}).

It then suffices to study the lower bound under the Frobenius norm. Similar
to the lower bound under the spectral norm one has
\begin{eqnarray*}
&&\inf_{\hat{\Sigma}}\max_{\theta\in\mathcal{F}^{\ast}}2_{\theta
}^{2}
\mathbb{E}_{%
\mathbf{X}|\theta} \bigl|\!\bigl|\!\bigl|\hat{\Sigma}-\Sigma( \theta)
\bigr|\!\bigr|\!\bigr|_{F}^{2}\\
&&\qquad\geq\min_{ \{ ( \theta,\theta
^{\prime
} ) \dvtx H ( \gamma( \theta) ,\gamma(
\theta
^{\prime} ) ) \geq1 \} }\frac{ |\!|\!|
\Sigma(
\theta) -\Sigma( \theta^{\prime} )
|\!|\!|_{F}^{2}}{%
H ( \gamma( \theta) ,\gamma( \theta^{\prime
} )
) }
\frac{p}{2}\min_{i}\llVert\bar\mathbb{P}_{i,0}
\wedge\bar\mathbb{P}_{i,1}\rrVert.
\end{eqnarray*}
It is easy to see%
\[
\min_{ \{ ( \theta,\theta^{\prime} ) \dvtx H (
\gamma(
\theta) ,\gamma( \theta^{\prime} ) ) \geq
1 \} }\frac{ |\!|\!|\Sigma( \theta) -\Sigma
(
\theta^{\prime} ) |\!|\!|_{F}^{2}}{H ( \gamma
( \theta
) ,\gamma( \theta^{\prime} ) ) }\asymp c_{n,p} \biggl(
\frac{\log p}{n} \biggr)^{1-q/2},
\]
and it follows from Lemma~\ref{affbd} that there is a constant $c>0$ such
that
\[
\min_{i}\llVert\bar\mathbb{P}_{i,0}\wedge\bar\mathbb{P}%
_{i,1}\rrVert\geq c. \]

\textit{Upper bound under Bregman matrix divergences}.
We now show that there exists an estimator $\hat{\Sigma}$ such that%
%
%
\begin{equation}
\mathbb{E}_{\mathbf{X}|\theta}\mathrm{L}_{\phi} ( \hat{\Sigma
},\Sigma)
\leq c \biggl[ c_{n,p} \biggl( \frac{\log p}{n} \biggr)^{1-q/2}+
\frac{1}{%
n} \biggr] \label{rateOper3}
\end{equation}
some constant $c>0$, uniformly over all $\phi\in\Phi$ and $\Sigma
\in
\mathcal{P}_{q}^{B}(\tau,c_{n,p})$. Let $A_{0}=\bigcap_{i,j}A_{ij}$,
where $%
A_{ij}$ is defined in (\ref{Aij}). Lemma~\ref{A0bound} yields that
%
%
\begin{equation}
\mathbb{P} ( A_{0} ) \geq1-2C_{1}p^{-5/2}.
\label{A0Tail}
\end{equation}
\upqed\end{pf}

%
\begin{lemma}
\label{BDTail}Let $\hat{\Sigma}_{B}$ be defined in equation (\ref%
{sigmahatStein}). Then for all $\Sigma\in\mathcal{G}_{q}^{B}(\rho
,c_{n,p}) $
\[
\mathbb{P} \biggl( \frac{\epsilon_{1}}{2}I\prec\hat{\Sigma
}_{B}\prec3
\tau I \biggr) \geq1-C_{1}p^{-5/2}.
\]
\end{lemma}

\begin{pf}
Write
$\hat{%
\Sigma}_{B}=\Sigma+ ( \hat{\Sigma}_{B}-\Sigma) $.
Since $%
|\!|\!|\hat{\Sigma}_{B}-\Sigma|\!|\!|\leq
|\!|\!|\hat{\Sigma}%
_{B}-\Sigma|\!|\!|_{1}$, the lemma is then a direct
consequence of
Lemma~\ref{A0bound} and equation (\ref{boundDominate}) which implies $
|\!|\!|\hat{\Sigma}_{B}-\Sigma|\!|\!|_{1}\leq
Cc_{n,p} ( \frac{%
\log p}{n} )^{ ( 1-q ) /2}\rightarrow0$ over
$A_{0}$. %
\end{pf}

Lemma~\ref{FroBD} implies%
%
%
\begin{eqnarray}\label{sparsel1bd}
\qquad &&\mathbb{E}_{\mathbf{X}|\theta}\mathrm{L}_{\phi} ( \hat{
\Sigma}%
_{B},\Sigma)\nonumber\\
&&\qquad=\mathbb{E}_{\mathbf{X}|\theta} \bigl\{
\mathrm{L}%
_{\phi} ( \hat{\Sigma}_{B},\Sigma)
I(A_{0}) \bigr\} +\mathbb{E}_{%
\mathbf{X}|\theta} \bigl\{
\mathrm{L}_{\phi} ( \hat{\Sigma}_{B},\Sigma) I
\bigl(A_{0}^{c}\bigr) \bigr\}
\nonumber
\\[-8pt]
\\[-8pt]
\nonumber
&&\qquad\leq C\mathbb{E}_{\mathbf{X}|\theta} \biggl\{ \frac{1}{p} |\!|\!|
\hat{%
\Sigma}-\Sigma|\!|\!|_{F}^{2}I(A_{0})
\biggr\} +\mathbb{E}_{\mathbf{X}%
|\theta} \bigl\{ \mathrm{L}_{\phi} ( \hat{
\Sigma}_{B},\Sigma) I\bigl(A_{0}^{c}\bigr) \bigr
\}
\\
&&\qquad\leq 16C\sup_{j}\sum_{i\neq j}\min\biggl\{
\vert\sigma_{ij}\vert^{2},\gamma\frac{\log p}{n}
\biggr\} +\mathbb{E}_{\mathbf{X}%
|\theta} \bigl\{ \mathrm{L}_{\phi} ( \hat{
\Sigma}_{B},\Sigma) I\bigl(A_{0}^{c}\bigr) \bigr
\} +C\frac{1}{n}.\nonumber
\end{eqnarray}
The second term in (\ref{sparsel1bd}) is negligible since%
\begin{eqnarray*}
\mathbb{E}_{\mathbf{X}|\theta}\mathrm{L}_{\phi} ( \hat{
\Sigma}%
_{B},\Sigma) \bigl\{ A_{0}^{c}
\bigr\} &\leq&C\cdot\bigl[ \max\{ \log n,\log p \} \bigr]^{\vert
r\vert}\cdot
\mathbb{P} \bigl( A_{0}^{c} \bigr)
\\
&\leq&C\cdot\bigl[ \max\{ \log n,\log p \} \bigr]^{\vert r\vert
}C_{1}p^{-5/4}\\
&=& o
\biggl( c_{n,p} \biggl( \frac{\log p}{%
n} \biggr)^{1-q/2}
\biggr)
\end{eqnarray*}
by applying the Cauchy--Schwarz inequality twice. We now consider the first
term in equation (\ref{sparsel1bd}). Set $k^{\ast}=\lfloor
c_{n,p} ( {%
\frac{n}{\log p}} )^{q/2}\rfloor$. Then we have%
\begin{eqnarray*}
\sum_{i\neq j}\min\biggl\{ \vert
\sigma_{ij}\vert^{2},\gamma\frac{\log p}{n} \biggr\} &
\leq&\gamma^{2} \biggl( \sum_{i\leq
k^{\ast
}}+\sum
_{i>k^{\ast}} \biggr) \min\biggl\{ \vert
\sigma_{
[ i ]
j}\vert^{2},\frac{\log p}{n} \biggr\}
\\
&\leq&C_{3}k^{\ast}\frac{\log p}{n}+C_{3}\sum
_{i>k^{\ast}} \biggl( \frac{%
c_{n,p}}{i}
\biggr)^{2/q}
\\
&\leq&C_{4} \biggl[ k^{\ast}\frac{\log p}{n}+c_{n,p}^{2/q}k^{\ast
}
\cdot\bigl( k^{\ast} \bigr)^{-2/q} \biggr] \\
&\leq&
C_{5}c_{n,p} \biggl( \frac{\log p}{%
n}
\biggr)^{1-q/2},
\end{eqnarray*}
which immediately yields equation (\ref{rateOper3}).

\begin{supplement}
\stitle{Supplement to ``Optimal rates of convergence for sparse
covariance matrix estimation''}
\slink[doi]{10.1214/12-AOS998SUPP} 
\sdatatype{.pdf}
\sfilename{aos998\_supp.pdf}
\sdescription{In this supplement we prove the additional technical
lemmas used in the proof of Lemma~\ref{affbd}.}
\end{supplement}

%
%

\printaddresses


\begin{thebibliography}{27}

\bibitem[\protect\citeauthoryear{Abramovich et~al.}{2006}]{Abretal06}
%
\begin{barticle}[mr]
\bauthor{\bsnm{Abramovich},~\bfnm{Felix}\binits{F.}},
\bauthor{\bsnm{Benjamini},~\bfnm{Yoav}\binits{Y.}},
\bauthor{\bsnm{Donoho},~\bfnm{David~L.}\binits{D.~L.}} \AND
\bauthor{\bsnm{Johnstone},~\bfnm{Iain~M.}\binits{I.~M.}}
(\byear{2006}).
\btitle{Adapting to unknown sparsity by controlling the false
discovery rate}.
\bjournal{Ann. Statist.}
\bvolume{34}
\bpages{584--653}.
\bid{doi={10.1214/009053606000000074}, issn={0090-5364}, mr={2281879}}
\bptok{imsref}%
\end{barticle}
%
\endbibitem

\bibitem[\protect\citeauthoryear{Assouad}{1983}]{Ass83}
%
\begin{barticle}[mr]
\bauthor{\bsnm{Assouad},~\bfnm{Patrice}\binits{P.}}
(\byear{1983}).
\btitle{Deux remarques sur l'estimation}.
\bjournal{C. R. Acad. Sci. Paris S\'er. I Math.}
\bvolume{296}
\bpages{1021--1024}.
\bid{issn={0249-6291}, mr={0777600}}
\bptok{imsref}%
\end{barticle}
%
\endbibitem

\bibitem[\protect\citeauthoryear{Bickel and Levina}{2008a}]{BicLev08N1}
%
\begin{barticle}[mr]
\bauthor{\bsnm{Bickel},~\bfnm{Peter~J.}\binits{P.~J.}} \AND
\bauthor{\bsnm{Levina},~\bfnm{Elizaveta}\binits{E.}}
(\byear{2008}a).
\btitle{Regularized estimation of large covariance matrices}.
\bjournal{Ann. Statist.}
\bvolume{36}
\bpages{199--227}.
\bid{doi={10.1214/009053607000000758}, issn={0090-5364}, mr={2387969}}
\bptok{imsref}%
\end{barticle}
%
\endbibitem

\bibitem[\protect\citeauthoryear{Bickel and Levina}{2008b}]{BicLev08N2}
%
\begin{barticle}[mr]
\bauthor{\bsnm{Bickel},~\bfnm{Peter~J.}\binits{P.~J.}} \AND
\bauthor{\bsnm{Levina},~\bfnm{Elizaveta}\binits{E.}}
(\byear{2008}b).
\btitle{Covariance regularization by thresholding}.
\bjournal{Ann. Statist.}
\bvolume{36}
\bpages{2577--2604}.
\bid{doi={10.1214/08-AOS600}, issn={0090-5364}, mr={2485008}}
\bptok{imsref}%
\end{barticle}
%
\endbibitem

\bibitem[\protect\citeauthoryear{Br{\`e}gman}{1967}]{Bre67}
%
\begin{barticle}[mr]
\bauthor{\bsnm{Br{\`e}gman},~\bfnm{L.~M.}\binits{L.~M.}}
(\byear{1967}).
\btitle{A relaxation method of finding a common point of convex sets
and its
application to the solution of problems in convex programming}.
\bjournal{USSR Comput. Math. Math. Phys.}
\bvolume{7}
\bpages{200--217}.
\bptok{imsref}%
\end{barticle}
%
\endbibitem

\bibitem[\protect\citeauthoryear{Cai and Liu}{2011}]{CaiLiu11}
%
\begin{barticle}[mr]
\bauthor{\bsnm{Cai},~\bfnm{Tony}\binits{T.}} \AND
\bauthor{\bsnm{Liu},~\bfnm{Weidong}\binits{W.}}
(\byear{2011}).
\btitle{Adaptive thresholding for sparse covariance matrix estimation}.
\bjournal{J. Amer. Statist. Assoc.}
\bvolume{106}
\bpages{672--684}.
\bid{doi={10.1198/jasa.2011.tm10560}, issn={0162-1459}, mr={2847949}}
\bptok{imsref}%
\end{barticle}
%
\endbibitem

\bibitem[\protect\citeauthoryear{Cai, Liu and Zhou}{2011}]{CaiLiuZho}
%
\begin{bmisc}[auto:STB|2012/10/26|14:52:12]
\bauthor{\bsnm{Cai},~\bfnm{T.~T.}\binits{T.~T.}},
\bauthor{\bsnm{Liu},~\bfnm{W.}\binits{W.}} \AND
\bauthor{\bsnm{Zhou},~\bfnm{H.~H.}\binits{H.~H.}}
(\byear{2011}).
\bhowpublished{Optimal estimation of large sparse precision matrices.
Unpublished manuscript}.
\bptok{imsref}%
\end{bmisc}
%
\endbibitem

\bibitem[\protect\citeauthoryear{Cai, Zhang and Zhou}{2010}]{CaiZhaZho10}
%
\begin{barticle}[mr]
\bauthor{\bsnm{Cai},~\bfnm{T.~Tony}\binits{T.~T.}},
\bauthor{\bsnm{Zhang},~\bfnm{Cun-Hui}\binits{C.-H.}} \AND
\bauthor{\bsnm{Zhou},~\bfnm{Harrison~H.}\binits{H.~H.}}
(\byear{2010}).
\btitle{Optimal rates of convergence for covariance matrix estimation}.
\bjournal{Ann. Statist.}
\bvolume{38}
\bpages{2118--2144}.
\bid{doi={10.1214/09-AOS752}, issn={0090-5364}, mr={2676885}}
\bptok{imsref}%
\end{barticle}
%
\endbibitem

\bibitem[\protect\citeauthoryear{Cai and Zhou}{2009}]{CaiZho}
%
\begin{bmisc}[auto:STB|2012/10/26|14:52:12]
\bauthor{\bsnm{Cai},~\bfnm{T.~T.}\binits{T.~T.}} \AND
\bauthor{\bsnm{Zhou},~\bfnm{H.~H.}\binits{H.~H.}}
(\byear{2009}).
\bhowpublished{Covariance matrix estimation under the $\ell_{1}$ norm (with
discussion). \textit{Statist. Sinica} \textbf{22} 1319--1378}.
\bptok{imsref}%
\end{bmisc}
%
\endbibitem

\bibitem[\protect\citeauthoryear{Cai and Zhou}{2012}]{supp}
%
\begin{bmisc}[auto]
\bauthor{\bsnm{Cai},~\bfnm{T.~T.}\binits{T.~T.}} \AND
\bauthor{\bsnm{Zhou},~\bfnm{H.~H.}\binits{H.~H.}}
(\byear{2012}).
\bhowpublished{Supplement to ``Optimal rates of convergence for sparse
covariance matrix estimation.'' DOI:\doiurl{10.1214/12-AOS998SUPP}}.
\bptok{imsref}%
\end{bmisc}
%
\endbibitem

\bibitem[\protect\citeauthoryear{Censor and Zenios}{1997}]{CenZen97}
%
\begin{bbook}[mr]
\bauthor{\bsnm{Censor},~\bfnm{Yair}\binits{Y.}} \AND
\bauthor{\bsnm{Zenios},~\bfnm{Stavros~A.}\binits{S.~A.}}
(\byear{1997}).
\btitle{Parallel Optimization: Theory, Algorithms, and Applications}.
\bpublisher{Oxford Univ. Press}, \blocation{New York}.
\bid{mr={1486040}}
\bptok{imsref}%
\end{bbook}
%
\endbibitem

\bibitem[\protect\citeauthoryear{Dhillon and Tropp}{2007}]{DhiTro07}
%
\begin{barticle}[mr]
\bauthor{\bsnm{Dhillon},~\bfnm{Inderjit~S.}\binits{I.~S.}} \AND
\bauthor{\bsnm{Tropp},~\bfnm{Joel~A.}\binits{J.~A.}}
(\byear{2007}).
\btitle{Matrix nearness problems with {B}regman divergences}.
\bjournal{SIAM J. Matrix Anal. Appl.}
\bvolume{29}
\bpages{1120--1146}.
\bid{doi={10.1137/060649021}, issn={0895-4798}, mr={2369287}}
\bptok{imsref}%
\end{barticle}
%
\endbibitem

\bibitem[\protect\citeauthoryear{Donoho and Liu}{1991}]{DonLiu91}
%
\begin{barticle}[mr]
\bauthor{\bsnm{Donoho},~\bfnm{David~L.}\binits{D.~L.}} \AND
\bauthor{\bsnm{Liu},~\bfnm{Richard~C.}\binits{R.~C.}}
(\byear{1991}).
\btitle{Geometrizing rates of convergence. {II}}.
\bjournal{Ann. Statist.}
\bvolume{19}
\bpages{633--667}.
\bid{doi={10.1214/aos/1176348114}, issn={0090-5364}}
\bptok{imsref}%
\end{barticle}
%
\endbibitem

\bibitem[\protect\citeauthoryear{El~Karoui}{2008}]{ElK08}
%
\begin{barticle}[mr]
\bauthor{\bsnm{El~Karoui},~\bfnm{Noureddine}\binits{N.}}
(\byear{2008}).
\btitle{Operator norm consistent estimation of large-dimensional sparse
covariance matrices}.
\bjournal{Ann. Statist.}
\bvolume{36}
\bpages{2717--2756}.
\bid{doi={10.1214/07-AOS559}, issn={0090-5364}, mr={2485011}}
\bptok{imsref}%
\end{barticle}
%
\endbibitem

\bibitem[\protect\citeauthoryear{Kulis, Sustik and
Dhillon}{2009}]{KulSusDhi09}
%
\begin{barticle}[mr]
\bauthor{\bsnm{Kulis},~\bfnm{Brian}\binits{B.}},
\bauthor{\bsnm{Sustik},~\bfnm{M{\'a}ty{\'a}s~A.}\binits{M.~A.}}
\AND
\bauthor{\bsnm{Dhillon},~\bfnm{Inderjit~S.}\binits{I.~S.}}
(\byear{2009}).
\btitle{Low-rank kernel learning with {B}regman matrix divergences}.
\bjournal{J. Mach. Learn. Res.}
\bvolume{10}
\bpages{341--376}.
\bid{issn={1532-4435}, mr={2485986}}
\bptok{imsref}%
\end{barticle}
%
\endbibitem

\bibitem[\protect\citeauthoryear{Lam and Fan}{2009}]{LamFan09}
%
\begin{barticle}[mr]
\bauthor{\bsnm{Lam},~\bfnm{Clifford}\binits{C.}} \AND
\bauthor{\bsnm{Fan},~\bfnm{Jianqing}\binits{J.}}
(\byear{2009}).
\btitle{Sparsistency and rates of convergence in large covariance matrix
estimation}.
\bjournal{Ann. Statist.}
\bvolume{37}
\bpages{4254--4278}.
\bid{doi={10.1214/09-AOS720}, issn={0090-5364}, mr={2572459}}
\bptok{imsref}%
\end{barticle}
%
\endbibitem

\bibitem[\protect\citeauthoryear{Le~Cam}{1973}]{LeCam73}
%
\begin{barticle}[mr]
\bauthor{\bsnm{Le~Cam},~\bfnm{L.}\binits{L.}}
(\byear{1973}).
\btitle{Convergence of estimates under dimensionality restrictions}.
\bjournal{Ann. Statist.}
\bvolume{1}
\bpages{38--53}.
\bid{issn={0090-5364}, mr={0334381}}
\bptok{imsref}%
\end{barticle}
%
\endbibitem

\bibitem[\protect\citeauthoryear{Le~Cam}{1986}]{LeCam86}
%
\begin{bbook}[mr]
\bauthor{\bsnm{Le~Cam},~\bfnm{Lucien}\binits{L.}}
(\byear{1986}).
\btitle{Asymptotic Methods in Statistical Decision Theory}.
\bpublisher{Springer}, \blocation{New York}.
\bid{doi={10.1007/978-1-4612-4946-7}, mr={0856411}}
\bptok{imsref}%
\end{bbook}
%
\endbibitem

\bibitem[\protect\citeauthoryear{Ravikumar et~al.}{2008}]{Ravetal08}
%
\begin{bmisc}[auto:STB|2012/10/26|14:52:12]
\bauthor{\bsnm{Ravikumar},~\bfnm{P.}\binits{P.}},
\bauthor{\bsnm{Wainwright},~\bfnm{M.}\binits{M.}},
\bauthor{\bsnm{Raskutti},~\bfnm{G.}\binits{G.}} \AND
\bauthor{\bsnm{Yu},~\bfnm{B.}\binits{B.}}
(\byear{2008}).
\bhowpublished{High-dimensional covariance estimation by minimizing
$l_{1}$-penalized log-determinant divergence. Technical Report 797, Dept.
Statistics, UC Berkeley.}
\bptok{imsref}%
\end{bmisc}
%
\endbibitem

\bibitem[\protect\citeauthoryear{Rothman, Levina and Zhu}{2009}]{RotLevZhu09}
%
\begin{barticle}[mr]
\bauthor{\bsnm{Rothman},~\bfnm{Adam~J.}\binits{A.~J.}},
\bauthor{\bsnm{Levina},~\bfnm{Elizaveta}\binits{E.}} \AND
\bauthor{\bsnm{Zhu},~\bfnm{Ji}\binits{J.}}
(\byear{2009}).
\btitle{Generalized thresholding of large covariance matrices}.
\bjournal{J. Amer. Statist. Assoc.}
\bvolume{104}
\bpages{177--186}.
\bid{doi={10.1198/jasa.2009.0101}, issn={0162-1459}, mr={2504372}}
\bptok{imsref}%
\end{barticle}
%
\endbibitem


\bibitem[\protect\citeauthoryear{Saulis and
Statulevi{\v{c}}ius}{1991}]{SauSta91}
%
\begin{bbook}[mr]
\bauthor{\bsnm{Saulis},~\bfnm{L.}\binits{L.}} \AND
\bauthor{\bsnm{Statulevi{\v{c}}ius},~\bfnm{V.~A.}\binits{V.~A.}}
(\byear{1991}).
\btitle{Limit Theorems for Large Deviations}.
\bseries{Mathematics and Its Applications (Soviet Series)}
\bvolume{73}.
\bpublisher{Kluwer Academic}, \blocation{Dordrecht}.
\bid{doi={10.1007/978-94-011-3530-6}, mr={1171883}}
\bptok{imsref}%
\end{bbook}
%
\endbibitem

\bibitem[\protect\citeauthoryear{Tsybakov}{2009}]{Tsybakov}
%
\begin{bbook}[mr]
\bauthor{\bsnm{Tsybakov},~\bfnm{Alexandre~B.}\binits{A.~B.}}
(\byear{2009}).
\btitle{Introduction to Nonparametric Estimation}.
\bpublisher{Springer}, \blocation{New York}.
\bid{doi={10.1007/b13794}, mr={2724359}}
\bptok{imsref}%
\end{bbook}
%
\endbibitem

\bibitem[\protect\citeauthoryear{van~der Vaart}{1998}]{van98}
%
\begin{bbook}[mr]
\bauthor{\bparticle{van~der} \bsnm{Vaart},~\bfnm{A.~W.}\binits{A.~W.}}
(\byear{1998}).
\btitle{Asymptotic Statistics}.
\bseries{Cambridge Series in Statistical and Probabilistic Mathematics}
\bvolume{3}.
\bpublisher{Cambridge Univ. Press}, \blocation{Cambridge}.
\bid{mr={1652247}}
\bptok{imsref}%
\end{bbook}
%
\endbibitem

\bibitem[\protect\citeauthoryear{Whittle}{1960}]{Whi60}
%
\begin{barticle}[mr]
\bauthor{\bsnm{Whittle},~\bfnm{P.}\binits{P.}}
(\byear{1960}).
\btitle{Bounds for the moments of linear and quadratic forms in independent
variables}.
\bjournal{Theory Probab. Appl.}
\bvolume{5}
\bpages{302--305}.
\bptok{imsref}%
\end{barticle}
%
\endbibitem

\bibitem[\protect\citeauthoryear{Yu}{1997}]{Yu97}
%
\begin{bincollection}[mr]
\bauthor{\bsnm{Yu},~\bfnm{Bin}\binits{B.}}
(\byear{1997}).
\btitle{Assouad, {F}ano, and {L}e {C}am}.
In \bbooktitle{Festschrift for {L}ucien {L}e {C}am}
(\beditor{\binits{D.} \bsnm{Pollard}},
\beditor{\binits{E.} \bsnm{Torgersen}}
\AND
\beditor{\binits{G.} \bsnm{Yang}}, eds.)
\bpages{423--435}.
\bpublisher{Springer}, \blocation{New York}.
\bid{mr={1462963}}
\bptok{imsref}%
\end{bincollection}
%
\endbibitem

\end{thebibliography}
\end{document}